\documentclass{article}
\usepackage{authblk}
\usepackage[utf8]{inputenc}
\usepackage[a4paper, total={7in, 9in}]{geometry}
\usepackage{subcaption}
\usepackage{graphicx}
\usepackage{multirow}

\usepackage{pgfplots}
\usepackage{amsmath,amsthm,bm}
\usepackage{amssymb}
\usepackage[toc,page]{appendix}
\usepackage{chngcntr}
\usepackage{caption}
\usepackage{subcaption}
\usepackage{multirow}
\usepackage{verbatim}
\usepackage[]{algorithm2e}
\usepackage{float}
\usepackage{tabularx}
\usepackage{graphicx}
\usepackage{indentfirst}
\usepackage{tikz-3dplot}

\newcommand{\tens}[1] {\bm{\mathsf{#1}}}
\newcommand{\tensc}[1] {{\mathsf{#1}}}
\newcommand{\btheta}{\bm{\theta}}
\title{Tensor train based isogeometric analysis for PDE approximation on parameter dependent geometries}
\author[1,2,*]{Ion Gabriel Ion}
\author[1,2]{Dimitrios Loukrezis}
\author[1,2]{Herbert De Gersem}
\affil[1]{Technische Universit\"at Darmstadt, Institute for Accelerator Science and Electromagnetic Fields (TEMF)}
\affil[2]{Technische Universit\"at Darmstadt, Centre for Computational Engineering}
\affil[*]{Corresponding author: ion@temf.tu-darmstadt.de}
\date{ }

\begin{document}
	
\maketitle

\begin{abstract}
\noindent This work develops a numerical solver based on the combination of isogeometric analysis (IGA) and the tensor train (TT) decomposition for the approximation of partial differential equations (PDEs) on parameter-dependent geometries.
First, the discrete Galerkin operator as well as the solution for a fixed geometry configuration are represented as tensors and the TT format is employed to reduce their computational complexity.
Parametric dependencies are included by considering the parameters that control the geometry configuration as additional dimensions next to the physical space coordinates.
The parameters are easily incorporated within the TT-IGA solution framework by introducing a tensor product basis expansion in the parameter space. 
The discrete Galerkin operators are accordingly extended to accommodate the parameter dependence, thus obtaining a single system that includes the parameter dependency.
The system is solved directly in the TT format and a low-rank representation of the parameter-dependent solution is obtained. 
The proposed TT-IGA solver is applied to several test cases which showcase its high computational efficiency and tremendous compression ratios achieved for representing the parameter-dependent IGA operators and solutions. \\

\noindent \textbf{keywords:} tensor decomposition, tensor train, isogeometric analysis, partial differential equations, parametric approximation, geometry deformation.
\end{abstract}
	
\section{Introduction}

Pioneered by Tom Hughes and collaborators, isogeometric analysis (IGA) is a powerful computational method for solving partial differential equations (PDEs) by combining finite element analysis (FEA) with computer aided design (CAD) \cite{cottrell2009isogeometric, hughes2005isogeometric, nguyen2015isogeometric}.
Using the IGA formulation, both the computational domain and the solution to the PDE are parametrized and represented by means of B-splines, non-uniform rational B-splines (NURBS), or tensor products thereof in more than one dimensions. 
Galerkin projection is then applied to derive a linear system of equations, the solution of which provides an approximate solution to the PDE.
A crucial advantage of IGA is the ability to represent complex geometries with minimal geometric approximation errors, due to the use of CAD tools.
Another important benefit to using IGA instead of the standard finite element method (FEM) is the comparatively reduced number of degrees of freedom (dofs) in relation to the accuracy of the approximate PDE solution \cite{da2014mathematical}. 
Owning to these advantages, IGA has been employed to resolve  problems governed by PDEs in numerous domains of application, exemplarily, structural analysis \cite{cottrell2006isogeometric, morganti2015patient, reali2006isogeometric, schmidt2010realization}, fluid mechanics \cite{akkerman2011isogeometric, bazilevs2008nurbs, garcia2019refined, hsu2011high, liu2019mixed, wang2017modeling}, and electromagnetics \cite{buffa2020isogeometric, buffa2010isogeometric,  dolz2019isogeometric, dolz2019numerical, simona2020isogeometric}.

On the downsides, a major disadvantage of IGA is the computationally expensive construction of the discrete Galerkin operators, in particular the mass and stiffness matrices. 
Especially challenging is the three-dimensional case, where the storage complexity is $\mathcal{O}(n^3p^3)$ in sparse format  \cite{hughes2005isogeometric}, where $p$ is the maximum degree of the individual B-spline functions and $n$ is the maximum size of the univariate B-spline bases. 
Moreover, the complexity of explicitly computing the matrix entries is bounded from below by the storage complexity. 
This bottleneck is greatly exacerbated in cases where IGA-based  models are utilized within parametric studies, such as shape optimization \cite{fusseder2015fundamental, merkel2021shape, pels2015optimization, wang2018structural}, uncertainty quantification (UQ) for stochastic geometry deformations \cite{georg2019uncertainty, zhang2019development}, or studies utilizing shape morphing techniques \cite{ziegler2022mode}. 
These parametric studies demand that multiple, often numerous, geometry configurations must be explored until an optimal shape or a statistical quantity of interest (QoI) can be estimated to sufficient accuracy. 
In turn, each geometry configuration corresponds to a computationally taxing re-assembly of the IGA system matrices, thus resulting to a possibly intractable computational cost for the full parametric study.

Several fast assembly procedures have been suggested in the literature, where the common idea is to exploit the Kronecker product \cite{van2000ubiquitous} (tensor product) structure of the IGA system matrices \cite{antolin2015efficient, hofreither2018black, mantzaflaris2015integration, mantzaflaris2017low, mantzaflaris2019low, pan2019low}. 
To reduce the cost of the multidimensional quadratures that are used to compute the stiffness and mass matrix entries, approaches based on sum factorization \cite{antolin2015efficient} and on integration by interpolation and look-up \cite{mantzaflaris2015integration} have been proposed.
The latter method was later extended with the use of partial tensor decompositions \cite{scholz2018partial}.
In \cite{mantzaflaris2017low} the stiffness and mass matrices are re-formatted as tensors and the canonical and Tucker tensor decompositions \cite{kolda2009tensor} are employed to reduce the computational complexity of the Galerkin matrix assembly.
A similar approach has been suggested for the parametrization of volumetric domains \cite{pan2019low}. 
An extension to space-time discretization of parabolic PDEs can be found in \cite{mantzaflaris2019low}, while the same idea is employed in \cite{bunger2020low} for PDE-constrained optimization, but using the so-called tensor train (TT) decomposition instead \cite{oseledets2011tensor, oseledets2009breaking}.
In \cite{hofreither2018black}, adaptive cross approximation (ACA) algorithms in two and three dimensions \cite{bebendorf2000approximation, bebendorf2011adaptive} are utilized to approximate the IGA stiffness matrix in a low-rank tensor format.
However, all aforementioned works consider physical computational domains with a fixed geometry.

This work focuses on the solution of PDEs on parameter-dependent geometries by means of IGA. 
In this context, different parameter realizations result in variations in the shape of the physical domain, hence, the full computational domain is a tensor space constructed by combining the physical and the parameter domains \cite{hackbusch2012tensor}. 
To address this challenging case, a framework is developed for incorporating parameter-dependent geometries in IGA solvers by exploiting the tensor-product structures and the corresponding multidimensional arrays, equivalently, tensors, that arise when discretizing the physical domain and the parameter space.
Parameter dependencies are accommodated by considering the geometrical parameters to be additional dimensions next to the physical coordinates of the problem under investigation. 
After both the physical domain and the parameter space have been suitably discretized, an approximation to the solution of the parameter-dependent PDE can be obtained by solving an extended system of equations, where the discrete operator (system matrix) as well as the discrete solution take the form of tensors.
To reduce the storage and computational complexity, these tensors are approximated and represented in the TT format, where further complexity reduction is accomplished by exploiting the particular structure of the discrete operator and using the so-called quantics TT (QTT) format \cite{khoromskij2011d, khoromskij2010quantics, oseledets2010approximation}. 
It is shown how the operator can be represented in the form of a TT-matrix and explicit expressions are given for the corresponding TT cores, i.e. the lower-dimensional tensors that the TT format consists of.
This step only requires the TT representation of the discretized, parameter-dependent geometry map, which is computed using TT-based approximation algorithms \cite{oseledets2010tt, savostyanov2011fast}.
Once the system operators are available in the TT format, a multilinear solver based on the alternating minimal energy (AMEn) algorithm \cite{dolgov2014alternating} is used to compute a low-rank approximation of the parameter-dependent PDE solution. 

While the TT-IGA solution framework proposed in this paper is novel, at least to the authors' knowledge, there exist previous works which have explored numerical methods for approximating PDEs on parameter-dependent geometries.
A rigorous study on the solution of PDEs on randomly deformed physical domains has been performed in \cite{castrillon2016analytic}, where a stochastic collocation method is put forth to approximate the PDE over the parameter space. 
Similar studies also based on stochastic collocation methods can be found in \cite{castrillon2021hybrid, harbrecht2008sparse, hiptmair2018large}.
In the same context, a TT-based stochastic Galerkin method is developed in \cite{eigel2020adaptive}, where however the discretization of the PDE is not based on IGA.
Similar to the latter work, the present paper also employs the TT format to reduce the complexity of the extended system matrices which are re-formatted as tensors, but now in the context of IGA-based discretizations. 
Additionally, the QTT format is employed to take advantage of the particular tensor-product structure of the spline-based geometry representations and further reduce the computational complexity.
Moreover, the framework developed in this work is not confined to random geometry deformations, but is more generally applicable, e.g. it can be applied for shape optimization or design space exploration without any modification.


The remaining of this paper is organized as follows. 
In section \ref{sec:premilinary}, tensor notation and the TT format are introduced, along with a brief presentation of the multilinear solver and the TT-based cross-approximation method employed in this work. 
A short presentation of B-splines is additionally available in the same section.
Next, in section \ref{sec:problem}, the IGA discretization is explained and the discrete operators of the Galerkin projection are given, along with the discretization of the parameter space.
In section \ref{sec:iga-tt}, the construction of the discrete IGA operators in the TT format is presented.  
Numerical results are presented and discussed in section \ref{sec:results}. 
First, a convergence study is performed on a geometry which is deformed based on a single parameter. 
The second study investigates the effectiveness of the proposed TT-IGA solver for an increasing number of parameters. 
In the third study, the TT-IGA solver is applied for a parameter-dependent domain with material jump discontinuity. 
In the final numerical study, the TT-IGA solver is applied for solving the scalar Helmholtz equation within a waveguide with varying geometrical features. 
Concluding remarks are available in section \ref{sec:conclusion}.


\section{Preliminaries and notation}
\label{sec:premilinary}
\subsection{Tensors and tensor notation}
In the context of this work, a \emph{tensor} is defined as a $d$-dimensional array and is denoted as $\tens{x} \in \mathbb{R}^{N_1 \times \cdots \times N_d}$. 
The dimensions of a tensor are also referred to as \emph{modes} or \emph{ways}, i.e. the expressions $d$-dimensional, $d$-way, and $d$-mode tensor are equivalent.
Tensor entries are denoted either as $\tensc{x}_{n_1 n_2 ... n_d}$, $n_k=1, \dots,N_k$, $k=1,\dots,d$, or using the multiindex notation $\tensc{x}_{\bm{n}}$, where $\bm{n}=(n_1, n_2,...,n_d)$ is a unique multiindex per tensor entry.
Tensor \emph{fibers} are the analog of matrix rows and columns and are denoted using the colon punctuation mark in the place of the corresponding index, e.g. $\tensc{x}_{n_1 n_2 \dots n_{k-1} : n_{k+1}, \dots, n_d}$ is the $k$-th fiber.
Tensor \emph{slices} of higher order are denoted in an analogous fashion.
A \emph{tensor-matrix}, also referred to as a \emph{tensor-operator}, is a generalization of the commonly used matrix-based operator to more than two dimensions and is denoted as $\tens{A} \in \mathbb{R}^{(M_1 \times \cdots \times M_d) \times (N_1 \times \cdots \times N_d)}$. 
The product between a tensor-matrix $\tens{A} \in \mathbb{R}^{(M_1 \times \cdots \times M_d) \times (N_1 \times \cdots \times N_d)}$ and a tensor $\tens{x} \in \mathbb{R}^{N_1 \times \cdots \times N_d}$ is defined similar to the standard matrix-vector product, and its result is a tensor of shape $M_1 \times \cdots \times M_d$, the entries of which are given by
\begin{equation}
\label{eq:tensor-product}
(\tens{Ax})_{\bm{m}} = \sum\limits_{\bm{n}} \tensc{A}_{\bm{m},\bm{n}} \tensc{x}_{\bm{n}}.
\end{equation}

\subsection{Tensor train decomposition}
The storage complexity of a tensor $\tens{x} \in \mathbb{R}^{N_1 \times \cdots \times N_d}$ is $\mathcal{O}\left(N^d\right)$, $N = \max_k N_k$, $k=1,\dots,d$, i.e. it scales exponentially with the tensor dimensions.
Tensor-based (multilinear) algebraic operations such as element-wise addition, multiplication, or summation over indices, scale with the same complexity.
To mitigate this so-called ``curse of dimensionality'', several tensor formats that reduce the complexity of tensor storage and tensor-based algebraic operations have been suggested in the literature \cite{kolda2009tensor}.

This work focuses on the so-called tensor train (TT) decomposition \cite{oseledets2011tensor}, where a $d$-dimensional tensor is represented using $d$ three-dimensional tensors. 
Using element-wise notation, a tensor $\tens{x}\in \mathbb{R}^{N_1\times\cdots\times N_d}$ is represented in the TT format as
\begin{equation}
\label{eq:tt-format}
\tensc{x}_{\bm{i}} = \sum\limits_{r_1=1}^{R_1}\sum\limits_{r_2=1}^{R_2}\cdots \sum\limits_{r_{d-1}=1}^{R_{d-1}} \tensc{g}^{(1)}_{1i_1r_1} \tensc{g}^{(2)}_{r_1i_2r_2} \cdots \tensc{g}^{(d)}_{r_{d-1}i_d1}, 
\end{equation}
where $\tens{g}^{(k)}\in\mathbb{R}^{R_{k-1} \times N_k \times R_k}$ are called the TT-cores and $\bm{R}=(1,R_1,...,R_{d-1},1)$ is the vector of the so-called TT-ranks. 
The obtained storage complexity is $\mathcal{O}(R^2Nd)$, i.e. the TT format provides linear complexity with respect to the tensor dimensions, as opposed to the exponential complexity of the full tensor.
Moreover, the complexity of basic multilinear algebraic operations scales linearly with the tensor dimension $d$ and polynomially with respect to the TT-ranks $\bm{R}$ and the mode sizes $N_1, \dots, N_d$ \cite{oseledets2011tensor}.

In general, an exact TT decomposition of a full tensor typically leads to high TT-ranks, thus increasing the computational complexity as well. 
However, in many cases, using a low-rank TT approximation $\tilde{\tens{x}}\approx \tens{x}$ is sufficient. 
Given the full tensor, a TT approximation can be computed using $d$ sequential singular value decompositions (SVDs) \cite{oseledets2011tensor, oseledets2009breaking}, up to a desired accuracy $\epsilon$ such that $\lVert\tens{x} - \tilde{\tens{x}}\rVert_{\text{F}} \leq \epsilon$, where $\lVert \cdot \rVert_{\text{F}}$ denotes the Frobenius norm. 
Moreover, reducing the rank of a TT-decomposition while maintaining a prescribed accuracy $\epsilon$, an operation called TT rounding, can be performed with the complexity $\mathcal{O}(R^3Nd)$.
Last but not least, assuming that the entries of a $d$-dimensional tensor stem from evaluating a $d$-dimensional function, the variables of which correspond to the tensor dimensions, e.g. $\tensc{x}_{n_1,\dots,n_d} = f(x_{1,n_1}, \dots, x_{d,n_d})$, $n_k=1,\dots,N_k$, $k=1,\dots,d$,  TT-based cross approximation algorithms can be employed to compute a TT approximation of the full tensor without ever constructing the full tensor \cite{oseledets2010tt, savostyanov2011fast}.

\subsection{Multilinear systems}
Of crucial interest in the context of this work is the solution of multilinear systems
\begin{equation}
\tens{A} \tens{x} = \tens{b},
\end{equation}
where $\tens{A} \in \mathbb{R}^{(N_1\times\cdots\times N_d)\times(N_1\times\cdots\times N_d)}$, $\tens{x}\in\mathbb{R}^{N_1\times\cdots\times N_d}$, and $\tens{b} \in \mathbb{R}^{N_1\times\cdots\times N_d}$.
In element-wise notation, the system is equivalently written as 
\begin{align}
	\sum_{\bm{n}} \tensc{A}_{\bm{m,n}} \tensc{x}_{\bm{n}} = \tensc{b}_{\bm{m}},
\end{align}
which obviously corresponds to the product between a tensor-matrix and a tensor defined in formula \eqref{eq:tensor-product}.
We assume that the operands $\tens{A}$ and $\tens{b}$ are given in the TT format and our goal is to represent the system's solution $\tens{x}$ in the TT-format as well. 
Generalizations of Krylov-subspace solvers have been proposed in the literature \cite{dolgov2013ttgmres}, however, a common bottleneck of these approaches is the typically large number of solver iterations which in turn lead to numerous computationally expensive rounding operations.
Alternatively, the solution can be computed by minimizing the residual of the system with respect to the TT-cores, which can be formulated as the nonlinear minimization problem
\begin{align}
\label{eq:tt-system-min}
\underset{\tens{g}^{(1)}, \dots, \tens{g}^{(d)}}{\min}||\tens{A}\tens{x}-\tens{b}||^2_\text{F}.
\end{align} 
However, if all TT-cores but one are fixed, the minimization \eqref{eq:tt-system-min} is transformed to a linear regression problem. 
The optimization can then be performed step-wise, such that only a single TT-core is optimized at each step, while all remaining cores remain fixed, a procedure known as alternating least squares (ALS) \cite{holtz2012alternating}.
The main disadvantage of the ALS algorithm is that the TT-ranks must be chosen a priori.
Alternatives that do not require a priori fixed TT-ranks have been proposed in the literature, such as the TT-based density matrix renormalization group (DMRG) \cite{oseledets2012solution} and the alternating minimal energy (AMEn) \cite{dolgov2014alternating} algorithms.
The DMRG algorithm proceeds similar to the ALS, but the residual minimization is performed for an extended TT-core $\tilde{\tens{g}}$, which is constructed by contracting two neighboring TT-cores such that the extended core's entries are given as
\begin{equation}
\tilde{\tensc{g}}_{r_{k-1}i_ki_{k+1}r_{k+1}} = \sum_{r_k}\tensc{g}^{(k)}_{r_{k-1}i_kr_{k}}\tensc{g}^{(k)}_{r_{k}i_{k+1}r_{k+1}}.
\end{equation} 
The optimized extended TT-core is then split into two separate cores using truncated SVD on a matricization of the extended core \cite{oseledets2012solution}. 
During the truncation step the rank can be adaptively chosen.
The AMEn algorithm is also similar to ALS, however, in each AMEn iteration a basis enlargement is performed in order to tackle the rank adaptivity issue. 
In comparison to the DMRG, the AMEn algorithm offers a more favorable complexity with respect to the tensor mode sizes \cite{dolgov2014alternating}.


\subsection{Computer aided design and B-splines}
\label{subsec:cad}
\begin{figure}[h]
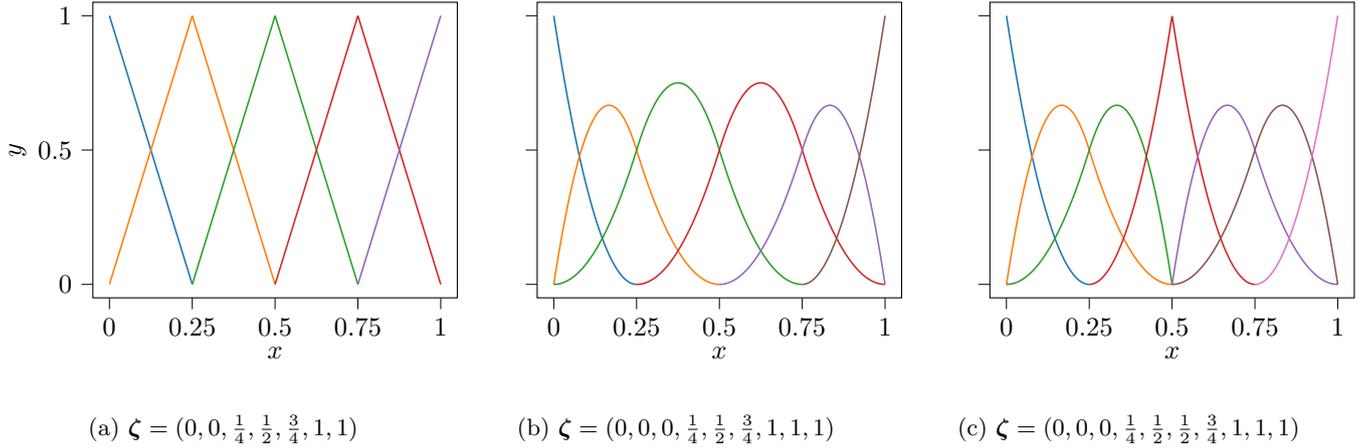

	\begin{subfigure}[b]{0.33\textwidth}
	\input{./images/bspl1}
		\caption{$\bm{\zeta}=(0,0,\frac{1}{4},\frac{1}{2},\frac{3}{4},1,1)$} \label{fig:bspl1}
	\end{subfigure}
	\begin{subfigure}[b]{0.33\textwidth}
	\input{./images/bspl2}
		\caption{$\bm{\zeta}=(0,0,0,\frac{1}{4},\frac{1}{2},\frac{3}{4},1,1,1)$} \label{fig:bspl2}
	\end{subfigure}
	\begin{subfigure}[b]{0.33\textwidth}
	\input{./images/bspl3}
		\caption{$\bm{\zeta}=(0,0,0,\frac{1}{4},\frac{1}{2},\frac{1}{2},\frac{3}{4},1,1,1)$} \label{fig:bspl20}
	\end{subfigure}
	\caption{Different B-spline bases for $p=1$, $p=2$, and $p=2$ with increased multiplicity for the knot $1/2$ (only the compact support is represented).}
	\label{fig:bspl}
\end{figure}

CAD-based geometry representations are commonly based on free-form curves which define a geometry via a projection map from a reference domain, commonly $\left[0,1\right]^d$, $d\in\left\{1,2,3\right\}$, to the physical domain. 
B-splines, as well as their generalization, NURBS, play a central role in this procedure \cite{piegl1996nurbs}.
Let $\bm{\zeta} = \left[\zeta_1, \dots, \zeta_{n+p+1}\right]$ with $\zeta_k\leq \zeta_{k+1}$, $\zeta_1=\cdots=\zeta_{p+1}=0$, and $\zeta_{n+1}=\cdots=\zeta_{n+p+1}=1$, be a knot vector, where $n$ denotes the space dimension and $p$ the polynomial degree.
Then, the B-spline basis functions $\{b_{k,p}\}_{k=1}^{n}$ of degree $p$ are defined via the Cox–de Boor recursion formula
\begin{align}
b_{k,p}(x) = \frac{x-\zeta_{k}}{\zeta_{k+1}-\zeta_{k}} b_{k,p-1}(x) + \frac{\zeta_{k+p+1}-x}{\zeta_{k+p+1}-\zeta_{k+1}} b_{k+1,p-1}(x), \quad p>0,
\end{align}
\begin{align}
b_{k,0}(x) = \begin{cases}
1, \quad x \in [\zeta_{k},\zeta_{k+1}) ,\\
0, \quad \text{otherwise}.
\end{cases}
\end{align} 
For the case $p>1$, the B-spline basis functions are piecewise polynomials of degree $p$ between the knots and $p-1$ times continuous differentiable at the knots (see Figure \ref{fig:bspl} for specific examples). 
However, if certain knots $\zeta_k \in (0,1)$ are repeated $p-1$ times, the space spanned by the B-splines contains functions that are not differentiable at those knots but only continuous (see Figure \ref{fig:bspl20}). 

Using B-splines as basis functions, any curve $C$ can be defined through a smooth map
\begin{equation}
F: \left[0,1\right] \rightarrow C.
\end{equation} 
Surfaces and volumes can be similarly defined, where in those cases the basis functions are constructed as B-spline tensor products (Kronecker products). 
For example, considering a three-dimensional geometry $\Omega \subset \mathbb{R}^3$ and given knot vectors $\bm{\zeta}_d$, degrees $p_d$, and number of basis functions $n_d$, $d=1,2,3$, the trivariate B-spline basis functions are given by the product
\begin{equation}
\label{eq:tp-b-splines}
b_{\bm{k}, \bm{p}}(\bm{x}) = b^{(1)}_{k_1, p_1}(x_1) b^{(2)}_{k_2, p_2}(x_2) b^{(3)}_{k_3, p_3}(x_3),
\end{equation}
where $\bm{k} = \left(k_1, k_2, k_3\right)$ and $\bm{p} = \left(p_1, p_2, p_3\right)$. 
Equivalently, the three-dimensional B-spline basis is obtained by the tensor product of the univariate B-spline bases, i.e.
\begin{equation}
\left\{b_{\bm{k}, \bm{p}}\right\} = \bigotimes_{d=1}^3\left\{b_{k_d, p_d}\right\},
\end{equation}
Accordingly, the three-dimensional knots, also referred to as the control points, form the tensor grid
\begin{equation}
\label{eq:tensorgrid}
\bm{Z} = \bigotimes_{d=1}^3 \bm{\zeta}_d.
\end{equation}
Volume parametrizations can then be constructed using functions from the space 
\begin{align}
\text{span}\{b_{\bm{k}, \bm{p}}(\bm{x})\: : \: k_1=1,...,n_1,k_2=1,...,n_2,k_3=1,...,n_3\},
\end{align}
and the corresponding map reads
\begin{equation}
G: \left[0,1\right]^3 \rightarrow \Omega.
\end{equation}
In the following, the index $p$, respectively, the multiindex $\bm{p}$, will be omitted for the simplicity of notation. That is, unless stated otherwise, a common degree $p$ is employed for all univariate B-spline bases.

\section{Model problem and discretization}
\label{sec:problem}
We consider a three-dimensional physical domain $\Omega(\bm{\theta}) \subset \mathbb{R}^3$, where $\bm{\theta} = (\theta_1, \dots, \theta_{N_{\text{p}}}) \in \Xi$ is a parameter vector that controls the shape of $\Omega$ and $\Xi = \Xi_1 \times \cdots \times \Xi_{N_{\text{p}}} \subset \mathbb{R}^{N_\text{p}}$ is a box domain which bounds the possible realizations of $\bm{\theta}$.
Assuming a spatial parametrization based on a tensor-product B-spline basis as defined in \eqref{eq:tp-b-splines}, the physical domain $\Omega(\bm{\theta})$ can be given via a map $G : [0,1]^3 \times \Xi \rightarrow \mathbb{R}^3$, such that 
\begin{align}
G_s(\bm{y},\bm{\theta}) \approx \sum\limits_{\bm{k}} \tensc{p}_{s\bm{k}}(\bm{\theta}) b_{\bm{k}}(\bm{y}), \label{eq:geom_approx}
\end{align}
where $G_s$, $s\in\{1,2,3\}$ are the dimensional components of the map $G$, $b_{\bm{k}}(\bm{y}) = b_{k_1}^{(1)}(y_1) b_{k_2}^{(2)}(y_2) b_{k_3}^{(3)}(y_3)$ is an element of the tensor-product B-spline basis, $\bm{y} \in \left[0,1\right]^3$, and
the tensor $\tens{p}\in\mathbb{R}^{3\times n_1\times n_2 \times n_3}$ contains the control points that form the tensor grid \eqref{eq:tensorgrid}. 
Accordingly, the tensor fiber $\tensc{p}_{:\bm{k}}$ represents one three-dimensional control point.
For every parameter $\bm{\theta} \in \Xi$,  $G(\cdot, \bm{\theta})$ is assumed to be a continuous and piecewise smooth function. 
Moreover, for every point $\bm{y}$ in the reference domain $[0,1]^3$, the function $G(\bm{y}, \cdot)$ is assumed to be sufficiently smooth for polynomial interpolation.

Next we consider the parameter-dependent boundary value problem (BVP)
\begin{subequations}
\label{eq:bvp}
\begin{align}
\nabla\cdot(\kappa(\cdot,\btheta) \nabla u(\cdot,\btheta)) +\rho u(\cdot,\btheta) = f(\cdot,\btheta),& \quad\text{in}\: \Omega(\btheta),\\
u(\cdot,\btheta) = g(\cdot,\btheta),& \quad\text{on}\: \Gamma_{\text{D}}(\btheta),\\
\partial_{\bm{\nu}} u(\cdot,\btheta) = 0 ,& \quad\text{on}\: \Gamma_{\text{N}}(\btheta),
\end{align}
\end{subequations}
where $\rho \in \mathbb{R}$, $\partial_{\bm{\nu}}$ denotes the normal derivative and $\Gamma_{\text{D}}(\btheta)$, $\Gamma_{\text{N}}(\btheta)$ are the Dirichlet and Neumann boundaries of $\Omega(\btheta)$, respectively. 
The coefficient function $\kappa \in L^\infty(\Omega(\btheta))$ is assumed to be uniformly bounded from below for every $\btheta \in \Xi$. Moreover, the right hand side belongs to $L^2$ for every parameter $\btheta$. 
Under these assumptions, the problem is well posed for every parameter $\btheta \in \Xi$ with the solution $u(\cdot, \btheta) \in H^1(\Omega(\btheta))$ \cite{Steinbach2007NumericalAM}.
 
\subsection{Spatial discretization via IGA}
\label{subsec:iga}
In IGA, the discrete solution space $V_{\bm{n},\btheta} \subset H^1({\Omega(\btheta)})$, $\bm{n} = \left(n_1, n_2, n_3\right)$, is constructed using the same B-spline basis employed for the geometry discretization, such that 
\begin{align}
V_{\bm{n},\btheta} = \text{span}\{b_{\bm{k}} \circ G^{-1}(\cdot,\btheta) \: : \: k_1 = 1,...,n_1\: ,\: k_2=1,...,n_2\: , \: k_3 =1,...,n_3\},
\end{align}
where $G^{-1}(\cdot,\btheta)$ denotes the inverse map for a fixed parameter $\btheta \in \Xi$. 
When discretizing the BVP \eqref{eq:bvp} using the Galerkin projection \cite{monk2003finite}, the following mass and stiffness matrices must be constructed:
\begin{align}
\tensc{M}_{\bm{m},\bm{k}}(\btheta) =a^\text{M}_{\btheta}(u_{\bm{m}}(\cdot,\btheta),v_{\bm{k}}(\cdot,\btheta))=\int\limits_{\Omega(\btheta)}    u_{\bm{m}}(\bm{x},\btheta)  v_{\bm{k}}(\bm{x},\btheta) \text{d} \bm{x},
\end{align}
\begin{align}
\tensc{S}_{\bm{m},\bm{k}}(\btheta) = a^\text{S}_{\btheta}(u_{\bm{m}}(\cdot,\btheta),v_{\bm{k}}(\cdot,\btheta))  = \int\limits_{\Omega(\btheta)}  \kappa(\bm{x},\btheta) \nabla u_{\bm{m}}(\bm{x},\btheta) \cdot \nabla v_{\bm{k}}(\bm{x},\btheta) \text{d} \bm{x},
\end{align}
where the trial and test functions $u_{\bm{m}}$ and $v_{\bm{k}}$, respectively, are the basis functions that span the solution space $V_{\bm{n},\btheta}$. 
Using the substitution rule for integrals, the parameter dependence in the integration domain can be moved inside the integral for both the stiffness and the mass matrices, such that
\begin{align}
\label{eq:mass-matrix}
\tensc{M}_{\bm{m},\bm{k}}(\btheta) = \int\limits_{[0,1]^3} b_{\bm{m}}(\bm{y}) b_{\bm{k}}(\bm{y}) \omega(\bm{y},\btheta) \text{d} \bm{y},\quad \omega(\bm{y},\btheta) = | \det D_{\bm{y}}G(\bm{y},\btheta)|,
\end{align}
\begin{align}
\label{eq:stiffness-matrix}
\tensc{S}_{\bm{m},\bm{k}}(\btheta) = \int\limits_{[0,1]^3} \nabla b_{\bm{m}}(\bm{y})^{\top} \bm{K}(\bm{y},\btheta) \nabla b_{\bm{k}}(\bm{y}) \hat{\kappa}(\bm{y},\btheta)\text{d}\bm{y},\quad \bm{K}(\bm{y},\btheta) =  D_{\bm{y}}G(\bm{y},\btheta)^{-\top} D_{\bm{y}}G(\bm{y},\btheta)^{-1} \omega(\bm{y},\btheta),
\end{align}
where $D_{\bm{y}}G(\bm{y},\btheta)$ is the Jacobian of the geometry mapping $G$ with respect to the reference coordinates and $\hat{\kappa}(\bm{y}, \btheta) = \kappa(G(\bm{y}, \btheta), \btheta)$ is the coefficient function represented in the reference domain. 
In the following, we denote with $\hat{u}(\cdot,\btheta) \in H^1([0,1]^3)$ the solution represented in the reference domain. The corresponding discrete space $\hat{V}_{\bm{n},\btheta}$ is obtained as the tensor product of the univariate B-spline bases.

\subsection{Parameter space discretization}
\label{subsec:param}
Thus far, the model problem is semidiscrete, i.e. no considerations regarding the discretization of the parameter space have been made.
To accommodate the parameter dependence, we collocate the geometry approximation from \eqref{eq:geom_approx} as well as the solution  of the BVP \eqref{eq:bvp} on a tensor-product grid of parameter realizations, $\Theta = \left\{\theta_{i_1}^{\text{C},1}\right\}_{i_1=1}^{\ell_1}\times\left\{\theta_{i_2}^{\text{C},2}\right\}_{i_2=1}^{\ell_2}\times\cdots\times\left\{\theta_{i_{N_{\text{p}}}}^{\text{C},N_{\text{p}}}\right\}_{i_{N_{\text{p}}}=1}^{\ell_{N_\text{p}}} \subset \Xi$, where the superscript ``$\text{C}$'' is used to denote the collocation points. 
In this work, the univariate grids consist of Gauss-Legendre nodes, however, other collocation points can be chosen as well, e.g. Chebyshev nodes. 
The geometry parametrization \eqref{eq:geom_approx} collocated on the grid ${\Theta}$ is modified to 
\begin{align}
	G_s(\bm{y},\bm{\theta}_{\bm{i}}) \approx \sum\limits_{\bm{k}} \tensc{p}_{s\bm{k}}(\bm{\theta}_{\bm{i}}) b_{\bm{k}}(\bm{y}) = \sum\limits_{\bm{k}} \tensc{p}_{s\bm{ki}} b_{\bm{k}}(\bm{y}), \label{eq:geom_approx2}
\end{align}
for $\btheta_{\bm{i}}=(\theta_{i_1}^{\text{C},1},...,\theta_{i_{N_{\text{p}}}}^{\text{C},N_{\text{p}}})$. 
The dofs of the fully discretized geometry map construct a $(1+3+N_{\text{p}})$-dimensional tensor $\tens{p}$ of size $3\times n_1\times n_2\times n_3\times \ell_1 \times\cdots\times \ell_{N_{\text{p}}}$, where the $4$-dimensional slices $\tensc{p}_{:\bm{i}}$ correspond to the dofs of the semidiscretized map for the parameter $\btheta_{\bm{i}}$, see formula \eqref{eq:geom_approx}. 
The fully discrete solution is similarly represented for every node $\btheta_{\bm{i}} \in \Theta$ as
\begin{align}
u(\bm{x},\btheta_{\bm{i}}) = \hat{u}(\bm{y},\btheta_{\bm{i}}) \approx \sum\limits_{\bm{k}}  \tensc{u}_{\bm{ki}} b_{\bm{k}}(\bm{y}).
\end{align}
The fully discrete solution takes the form of a tensor $\tens{u}\in\mathbb{R}^{n_1\times n_2\times n_3\times \ell_1 \times\cdots\times \ell_{N_{\text{p}}}}$ and the continuous approximation of the solution can be (approximately) recovered using interpolation, such that
\begin{align}
 \hat{u}(\bm{y},\btheta) \approx \hat{u}_{\bm{n},\bm{\ell}}(\bm{y},\btheta) = \sum\limits_{\bm{k},\bm{i}}  \tensc{u}_{\bm{ki}} b_{\bm{k}}(\bm{y}) P^{(1)}_{i_1}(\theta_1)\cdots P^{(N_{\text{p}})}_{i_{N_{\text{p}}}}(\theta_{N_{\text{p}}})= \sum\limits_{\bm{k},\bm{i}}  \tensc{u}_{\bm{ki}} b_{\bm{k}}(\bm{y}) P_{\bm{i}}(\btheta),
\end{align}
where $\hat{u}_{\bm{n},\bm{\ell}}$ denotes the approximate fully discrete solution and $\left\{P_{i_k}^{(k)}\right\}_{i_k=1}^{\ell_k}$ are univariate Legendre bases corresponding to the collocation points $\left\{\theta_{i_k}^{\text{C},k}\right\}_{i_k=1}^{\ell_k}$. 

Formally, the discrete solution $\hat{u}$ belongs to the tensor-product space $\hat{V}_{\bm{n}} \otimes (\mathcal{P}_{\ell_1} \otimes \cdots \otimes \mathcal{P}_{\ell_{N_{\text{p}}}})$, where $\hat{V}_{\bm{n}}$ is the space spanned by the multivariate B-splines and $\mathcal{P}_{\ell_{k}}$ are spaces of Legendre polynomials up to degree $\ell_{k}-1$.
The recovery of the solution can be achieved by solving $\prod_{k=1}^{N_{\text{p}}} \ell_k$ multilinear systems
\begin{align}
\sum\limits_{\bm{k}} \tensc{L}_{\bm{m},\bm{k}}(\btheta_{\bm{i}}) \tensc{u}_{\bm{ki}} = \tensc{f}_{\bm{m}}(\btheta_{\bm{i}}), \quad \forall \bm{m},\bm{i},
\end{align}
where $\tens{L}(\btheta)\in\mathbb{R}^{(n_1\times n_2 \times n_3) \times (n_1 \times n_2 \times n_3)}$ is the discrete operator obtained from the Galerkin discretization and enforcing the boundary conditions, i.e. it includes both the stiffness and the mass terms.  
Alternatively, a single extended system can be derived by constructing an $(n_1\times n_2\times n_3 \times \ell_1 \times \cdots \times \ell_{N_{\text{p}}})\times (n_1\times n_2\times n_3 \times \ell_1 \times \cdots \times \ell_{N_{\text{p}}})$ tensor-matrix and the corresponding right hand side, such that
\begin{align}
\tensc{L}_{\bm{mi},\bm{kq}} = \delta_{\bm{i}}^{\bm{q}}\tensc{L}_{\bm{m},\bm{k}}(\btheta_{\bm{i}}), \quad \tensc{f}_{\bm{mi}} = \tensc{f}_{\bm{m}}(\btheta_{\bm{i}}),
\end{align}
where $\delta_{\bm{i}}^{\bm{q}}$ denotes the (multidimensional) Kronecker delta. Finally, a multilinear system is solved to obtain the fully discrete solution $\tens{u} \in \mathbb{R}^{(n_1\times n_2\times n_3) \times (\ell_1 \times \cdots \times \ell_{N_{\text{p}}})}$, such that 
\begin{align}
\sum\limits_{\bm{k,q}} \tensc{L}_{\bm{mi},\bm{kq}} \tensc{u}_{\bm{kq}} = \tensc{f}_{\bm{mi}}, \quad \forall \bm{m},\bm{i}. \label{eq:multilinear_system}
\end{align}

We introduce the following space for parameter-dependent functions defined on the physical domain
\begin{align}
L^2(\Xi,G,\alpha) = \{ u(\cdot,\cdot) \: : \: u(\cdot,\btheta) \in H^\alpha(\Omega(\btheta)) \; \forall \btheta \in \Xi\text{ and } ||u||^2_{L^2(\Xi,G,\alpha)}<\infty\},\quad \alpha \in \mathbb{Z},
\end{align}
with the corresponding norm
\begin{align}
||u||^2_{L^2(\Xi,G,\alpha)}=\int\limits_{\Xi} || u ||_{H^\alpha(\Omega(\btheta))}^2 \text{d} \bm{\theta}. \label{eq:norm_def}
\end{align}
Then, the solution represented in the reference domain belongs to a tensor-product space $H^1([0,1]^3)\otimes L^2(\Xi)$.
In \cite{castrillon2016analytic}, a priori error estimates are derived for the case of random domain deformations combined with the FEM for elliptic PDEs. 
The main idea consists of using a map between the reference domain and the physical domain. 
This map includes the dependence of the outcome on the randomly deformed geometry. 
This parametric dependence is moved inside a density term in the integrals arising from the FEM discretization, thus reducing the problem to the particular case studied in \cite{babuska2007stochastic}. 
The parameter-dependent IGA setting presented in this work follows exactly the same path due to the parameter-dependent geometry map $G$ defined in \eqref{eq:geom_approx}, respectively in \eqref{eq:geom_approx2}. 
The error between the actual solution and the approximation in the norm defined in \eqref{eq:norm_def} can be split as
\begin{align}
||\hat{u}-\hat{u}_{\bm{n,\ell}}||_{L^2(\Xi,G,\alpha)} \leq \epsilon_{\text{IGA}}(\bm{n})+\epsilon_{\Xi}(\bm{\ell}),
\end{align}
where $\epsilon_{\text{IGA}}$ is the error arising from the IGA discretization and $\epsilon_{\Xi}$ is the error arising from the parameter space discretization \cite{castrillon2016analytic}. 
The error term $\epsilon_{\text{IGA}}$ decreases with  $\mathcal{O}(n^{-(p+1)})$ for uniformly spaced B-spline bases of degree $p$ \cite{da2014mathematical}, while it holds that $\epsilon_{\Xi}\leq C \sum_k e^{-\alpha _k\ell_k}$ for positive $\alpha_k$ and $C>0$ \cite{babuska2007stochastic, castrillon2016analytic}.


\section{Low-rank representation of IGA operators}
\label{sec:iga-tt}
In section \ref{sec:problem}, the problem setup and the parameter discretization were introduced, where the discrete solution as well as the discrete operators are represented in tensor format, taking advantage of the tensor-product structure of the solution space. 
However, the additional dimensions due to accommodating the parameter dependence, drastically increase the storage and computational requirements for both the solution and the operators. 
The TT decomposition is therefore employed to render the representation of the tensors affordable. 
To that end, an efficient way for constructing the discrete stiffness and mass tensor-matrices $\tens{S}$ and $\tens{M}$, respectively, as well as the right hand side $\tens{f}$ of the multilinear system \eqref{eq:multilinear_system} is presented in the following. 
Once all tensors are represented in the TT format, a multilinear solver is employed to recover a TT approximation of the fully discrete parameter-dependent solution $\tens{u}$.

\subsection{Geometry interpolation}
\label{subsec:geometry-interpolation}
Before computing the discrete operators in the TT format, the tensors corresponding to the geometry representations given in  \eqref{eq:geom_approx2} must first be represented in the TT format.
Using the cross-approximation method \cite{oseledets2010tt}, the three components of the map $G$, here denoted with $s\in\{1,2,3\}$, are evaluated on the joint Greville-parameter grid \cite{johnson2005higher}, thus resulting in the corresponding tensors $\tens{g}^s$ with entries
\begin{align}
\tensc{g}^s_{\bm{mi}} = G_s\left(y^{\text{G},1}_{m_1},y^{\text{G},2}_{m_2},y^{\text{G},3}_{m_3},\theta_{i_1}^{\text{C},1},...,\theta_{i_{N_{\text{p}}}}^{\text{C},{N_{\text{p}}}}\right), \label{eq:greville1}
\end{align}
where $\{y^{\text{G},1}_{m_1}\}_{m_1}$, $\{y^{\text{G},2}_{m_2}\}_{m_2}$, $\{y^{\text{G},3}_{m_3}\}_{m_3}$ are the univariate Greville abscissae corresponding to the B-spline bases \cite{johnson2005higher}. 
The Greville abscissae corresponding to a B-spline basis $\{b_k\}_k$ with basis functions of degree $p$ and the knot vector $\bm{\zeta}$ are defined as $y^{\text{G}}_m= (\zeta_{m+1}+\cdots +\zeta_{m+p+1})/p$ and have the property that the matrix with entries $(b_k(y^{\text{G}}_m))_{km}$ is nonsingular.
The following systems can be then solved to obtain the slices $\tens{p}_{s::}$, $s \in \{1,2,3\}$, of the control points tensor:
\begin{align}
	\tens{B} \tens{p}_{s::} = \tens{g}^s, \quad \tens{B}_{\bm{mi},\bm{kq}}= b_{k_1}(y^{\text{G},1}_{m_1}) b_{k_2}(y^{\text{G},2}_{m_2}) b_{k_3}(y^{\text{G},3}_{m_3}) \delta_{\bm{i}}^{\bm{q}}.
\end{align}
A similar procedure can be applied to interpolate a function $f(\cdot,\btheta)\in H^1(\Omega(\btheta))$ by evaluating it on the Greville-parameter grid defined in \eqref{eq:greville1}, which can be accomplished using an adaptive cross approximation method \cite{dolgov2019parametric, oseledets2010tt, savostyanov2011fast}. 
The following system is then solved to obtain the right hand side $\tens{f}$ of the discrete representation:
\begin{align}
\sum\limits_{\bm{n,q}}\tensc{B}_{\bm{mi},\bm{nq}}\tensc{f}_{\bm{nq}} = f(\tensc{g}^1_{\bm{mi}},\tensc{g}^2_{\bm{mi}},\tensc{g}^3_{\bm{mi}}).
\end{align}
Since the TT-operator $\tens{B}$ has TT-rank $\bm{R}=\bm{1}$, solving this system is computationally inexpensive.

\subsection{Discrete operators}
We first consider the construction of the mass tensor-matrix. The integration over $[0,1]^3$ is performed by constructing a tensor-product grid of univariate quadrature points $\left\{y_{i_1}^{\text{Q},1}\right\}_{i_1}^{}\times\left\{y_{i_2}^{\text{Q},2}\right\}_{i_2}^{}\times\left\{y_{i_3}^{\text{Q},3}\right\}_{i_3}^{}$ and the corresponding quadrature weights $\tensc{w}_{\bm{i}}=w^{(1)}_{i_1}  w^{(2)}_{i_2}w^{(3)}_{i_3}$. 
The univariate quadrature grids are chosen according to the corresponding B-spline basis, such that they allow for exact polynomial integration in between the knots of the basis. 
The mass tensor-matrix is then approximated as
\begin{align}
	\tensc{M}_{\bm{mi},\bm{kq}} = \delta_{\bm{q}}^{\bm{i}} \int\limits_{U}\omega(\bm{y},\bm{\theta}^{\text{C}}_{\bm{i}}) b_{\bm{m}}(\bm{y}) b_{\bm{k}}(\bm{y}) \text{d}\bm{y } \approx \delta_{\bm{q}}^{\bm{i}} \sum\limits_{\bm{j}} \tensc{w}_{\bm{j}} \tensc{o}_{\bm{ji}}b_{{m}_1}({y}^{\text{Q},1}_{j_1})  \cdots b_{{m}_3}({y}^{\text{Q},3}_{j_3})b_{{k}_1}({y}^{\text{Q},1}_{j_1})  \cdots b_{{k}_3}({y}^{\text{Q},3}_{j_3}) ,
	\end{align}
where the tensor $\tens{o}$ contains the value of the function $\omega$ defined in formula \eqref{eq:mass-matrix}, evaluated on the Cartesian product between the quadrature grid and the collocation grid. 
If the control points of the geometry discretization are given in the TT format as described in section~\ref{subsec:geometry-interpolation}, then the tensor $\tens{o}$ can also be represented in the TT format as well, such that
\begin{align}
	\tensc{o}_{\bm{ji}} = \omega\left(y_{j_1}^{\text{Q},1},y_{j_2}^{\text{Q},2},y_{j_3}^{\text{Q},3},\btheta_{\bm{i}}\right) = \sum\limits_{\text{perm. }\sigma } \text{sgn}(\sigma)\prod\limits_{s=1}^{3} \left(D_{\bm{y}}G(\bm{y}^{\text{Q}}_{\bm{j}},\btheta_{\bm{i}})\right)_{s \sigma(s)}=\sum\limits_{\text{perm. }\sigma } \text{sgn}(\sigma)\prod\limits_{s=1}^{3} \left((\partial_{\sigma(s)}\tens{B}) \tens{p}_{s:\bm{i}}\right)_{\bm{j}},
\end{align}
where $\sigma(s)$ are permutations of the tuple $(1,2,3)$ and $\left(\partial_{\sigma(s)}\tens{B})_{\bm{j},\bm{i}} = \partial y_{\sigma(s)} b_{\bm{i}}(y_{j_1}^{\text{Q},1},y_{j_2}^{\text{Q},2},y_{j_3}^{\text{Q},3}\right) $ are the components of a rank-$\bm{1}$ TT-operator. 
Finally, the TT decomposition of the discrete mass operator is given by
\begin{align}
\tensc{M}_{\bm{mi},\bm{kq}} = \sum\limits_{\bm{r}} \prod\limits_{s=1}^{3}  \underbrace{ \left(\sum\limits_{j_s}\tensc{o}^{(s)}_{r_{s-1} j_s r_s}{w}^{(s)}_{j_s}b_{{m}_s}({y}^{\text{Q},s}_{j_s})b_{{k}_s}({y}^{\text{Q},s}_{j_s})\right)}_{M^{(s\leq 3)}} \left(\prod\limits_{s=4}^{3+N_{\text{p}}}\underbrace{\tensc{o}^{(s)}_{r_{s-1} i_{s-3} r_{s}} \delta_{i_{s-3}}^{q_{s-3}}}_{M^{(s>3)}}\right),
\end{align}
where $\tens{o}^{(s)}$ are the TT cores of the tensor $\tens{o}$.
The stiffness tensor-matrix can be constructed in the TT format in a similar way, such that
\begin{align}
\label{eq:stifness-tensor}
\tensc{S}_{\bm{mi},\bm{kq}} &= \delta_{\bm{q}}^{\bm{i}}\int\limits_{[0,1]^3} \nabla b_{\bm{m}}(\bm{y})^{\top} \bm{K}(\bm{y},\bm{\theta}^{\text{C}}_{\bm{i}}) \nabla b_{\bm{k}}(\bm{y}) \hat{\kappa}(\bm{y},\bm{\theta}^{\text{C}}_{\bm{i}})\text{d}\bm{y} \nonumber \\
&=  \delta_{\bm{q}}^{\bm{i}} \sum\limits_{\alpha,\beta=1}^3 \sum\limits_{\bm{j}} \tensc{w}_{\bm{j}} \partial_{y_{\alpha}} b_{\bm{m}}({y}^{\text{Q},1}_{j_1},{y}^{\text{Q},2}_{j_2},{y}^{\text{Q},3}_{j_3}) 
\partial_{y_{\beta}} b_{\bm{k}}({y}^{\text{Q},1}_{j_1},{y}^{\text{Q},2}_{j_2},{y}^{\text{Q},3}_{j_3}) \tensc{K}^{(\alpha, \beta)}_{\bm{ji}} \tensc{k}_{\bm{ji}},
\end{align}
where $\tensc{K}_{\bm{ji}}^{(\alpha, \beta)}=K_{\alpha\beta}({y}^{\text{Q},1}_{j_1},{y}^{\text{Q},2}_{j_2},{y}^{\text{Q},3}_{j_3},\theta_{i_1}^{\text{C},1}.,...,\theta_{i_{N_{\text{p}}}}^{\text{C},N_{\text{p}}})$ and $\tensc{k}_{\bm{ji}}=\hat{\kappa}({y}^{\text{Q},1}_{j_1},{y}^{\text{Q},2}_{j_2},{y}^{\text{Q},3}_{j_3},\theta_{i_1}^{\text{C},1}.,...,\theta_{i_{N_{\text{p}}}}^{\text{C},N_{\text{p}}})$. 
Note that the construction of the latter tensor in the TT format can be performed using a TT-based cross approximation method \cite{dolgov2019parametric}.
Compared to the mass matrix, assembling the stiffness matrix includes the elementwise inversion of the tensor $\tens{o}$, which is performed using the AMEn algorithm \cite{dolgov2014alternating}. 
Moreover, summing over the components of the matrix $\bm{K}$ increases the TT rank of the stiffness tensor $\tens{S}$, thus leading to a higher assembly time. 
The size of the cores scales with $\mathcal{O}(n^2R^2)$, however, as it turns out from the construction process, the TT cores have a band diagonal structure with respect to the inner two modes. During the construction process, we benefit from this structure by storing only the nonzero elements. 
However, the TT solver still requires the full format of the cores.

\subsection{Quantized tensor train (QTT) decomposition}
One way to speed up the computations in the TT format is to use the so-called quantized TT (QTT) decomposition \cite{khoromskij2011d}. 
The basic idea of the QTT format is to reshape a given tensor into a higher-dimensional one while simultaneously reducing the mode size. 
Let $\tens{x}\in \mathbb{R}^{N_1 \times \cdots \times N_d}$ be a tensor with mode sizes $N_k$, $k=1,...,d$, such that $\log_2 N_k \in \mathbb{N}$. 
The tensor $\tens{x}$ can then be reshaped into a $\left(\sum_k\log_2 N_k\right)$-dimensional tensor and then be represented in the TT format. 
In many cases, this prior transformation in terms of tensor dimensions leads to a better storage and computational complexity compared to applying the TT decomposition to the original tensor \cite{khoromskij2011d, khoromskij2010quantics, ion2021tensor}. 
If the tensor modes are not powers of 2, the reshaping operation can still be performed using the prime number decomposition of the individual modes. 
In the context of this work, the QTT format is used to speed up the construction of the stiffness tensor-matrix, in particular regarding the elementwise inversion of the tensor $\tens{o}$ which is necessary for computing the tensors $\tens{K}^{(\alpha,\beta)}$, see formula \eqref{eq:stifness-tensor}. 
Additionally, the QTT format can be used when solving the multilinear system \eqref{eq:multilinear_system}, which often results in computational gains as well.  

\section{Numerical Results}
\label{sec:results}
In the following numerical investigations, we employ the TT-IGA framework proposed in this paper to solve PDEs on parameter-dependent geometries.  Comparisons against alternative solution methods are also performed. 
The first numerical investigation presents a convergence analysis for the case of single parameter dependence. 
The second example showcases how the TT-IGA solver scales with an increasing number of parameters.
The third test case concerns a parameter-dependent material jump within the computational domain. 
In the final numerical example, the scalar Helmholtz equation is solved within a waveguide structure subject to geometry deformations.
All computations have been performed on a standard workstation with a 20-core/40-thread Intel Xeon CPU, 2.2 GHz, and 96GB RAM. 
The \texttt{Python} package \texttt{torchtt}\footnote{\texttt{https://github.com/ion-g-ion/torchTT}}, which has been developed as part of this work, was used for all TT-based multilinear algebraic operations and system solutions.

\subsection{Test case 1: TT-IGA solver performance for single parameter dependence}
\label{subsec:test-case-1}

\begin{figure}[t!]
	\centering
	\begin{subfigure}[t]{0.49\textwidth}
		\includegraphics[width=\textwidth]{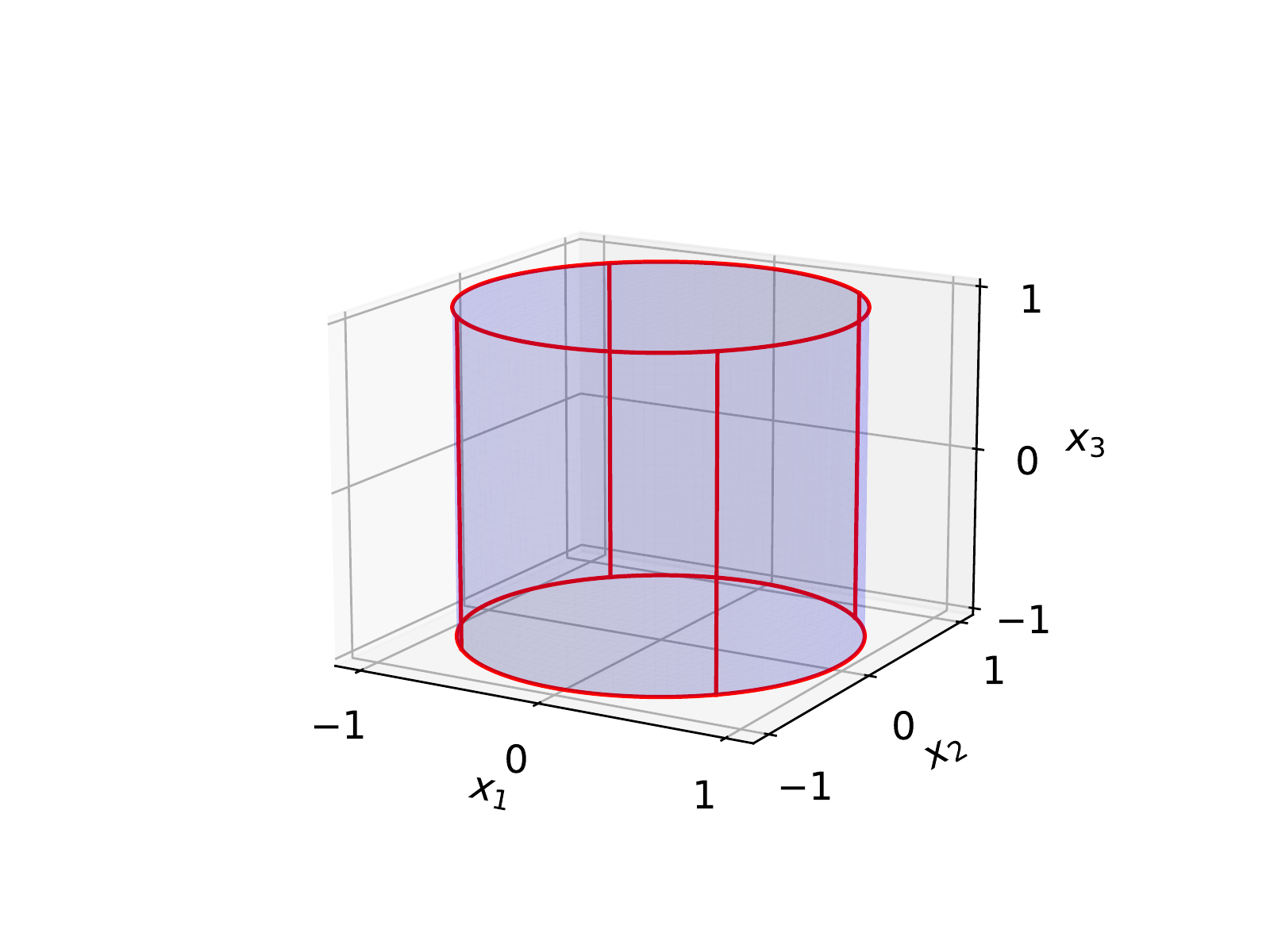}
		\caption{Nominal geometry $\left(\theta=0\right)$.}
	\end{subfigure}
	\begin{subfigure}[t]{0.49\textwidth}
		\includegraphics[width=\textwidth]{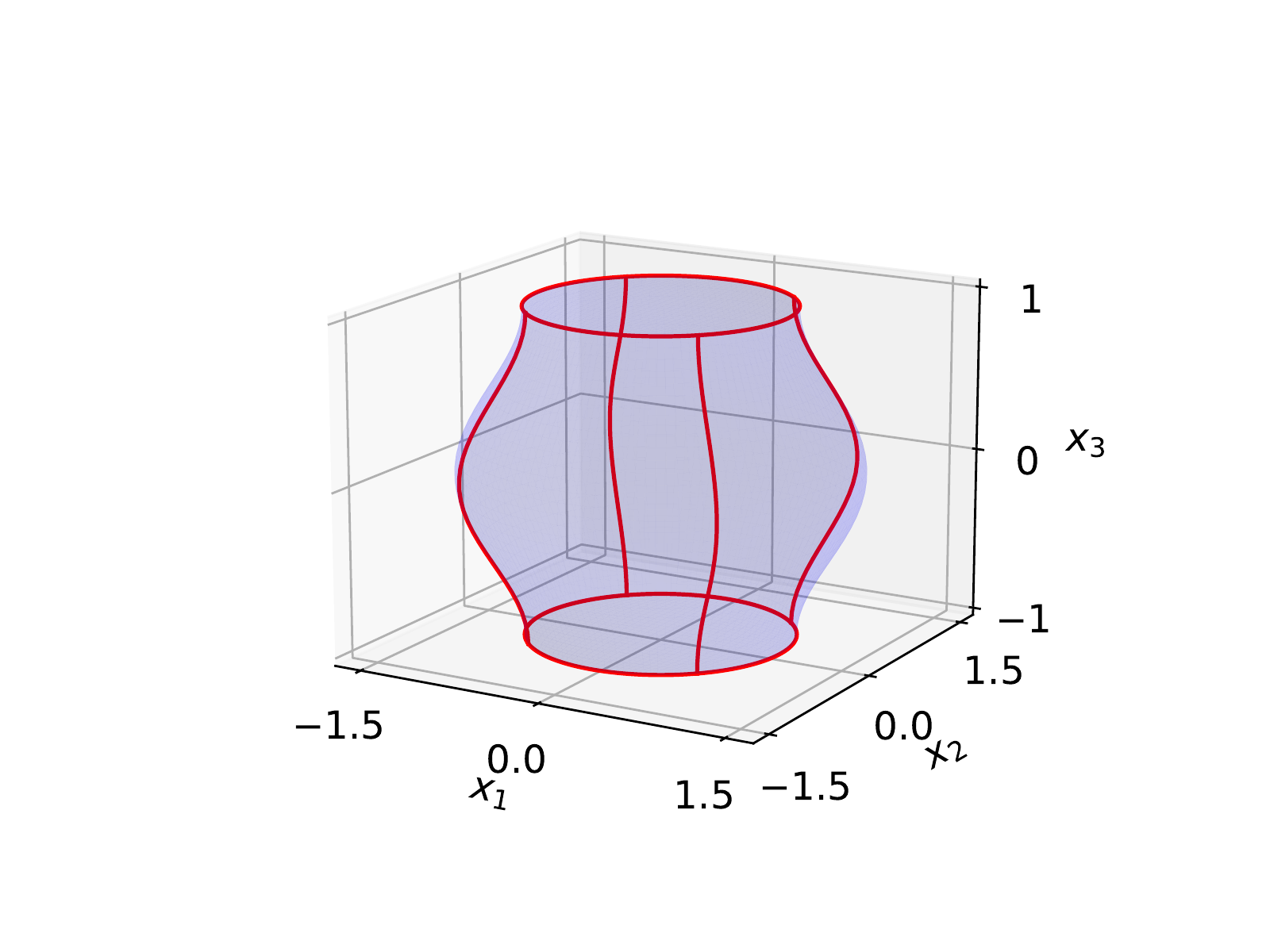}
		\caption{Fully deformed geometry $\left(\theta=1\right)$.}
	\end{subfigure}
	\caption{Test case 1: Deformation of a cylinder depending on a single parameter $\theta \in \left[0,1\right]$.}
	\label{fig:cyl2cav}
\end{figure}

We first perform a numerical study to assess the performance of the proposed TT-IGA solution method in terms of convergence and computational cost. 
The latter is separated into simulation runtime, storage needs, and time needed for assembling the stiffness matrix.
For this study, the computational domain is a cylinder which is deformed by varying a single parameter. 
Figure~\ref{fig:cyl2cav} shows the initial and the fully deformed geometry, which is given by the parametrization 
\begin{align}
	G(\bm{y},\theta)=\begin{pmatrix}
		(2y_1-1)\sqrt{1-\frac{(2y_2-1)^2}{2}}\left(\frac{1+\cos((2y_3-1)\pi)}{2}\theta+1\right) \\
		(2y_2-1)\sqrt{1-\frac{(2y_1-1)^2}{2}}\left(\frac{1+\cos((2y_3-1)\pi)}{2}\theta+1\right) \\
		2y_3-1
	\end{pmatrix}, \quad y_1,y_2,y_3,\theta \in [0,1].
\end{align}
In the parameter-dependent cylindrical domain $\Omega(\theta)$, the Poisson equation 
\begin{subequations}
\begin{align}
	\Delta u(\bm{x},\theta) &= 0 ,  &&\bm{x}\in\Omega(\theta), \\
	u(\bm{x},\theta) &= \sin(2x_1) \sin(3x_2) \exp(-\sqrt{13}x_3),   &&\bm{x}\in\partial\Omega(\theta),
\end{align}
\end{subequations}
is solved, the analytical solution of which is $u(\bm{x},\theta) = \sin(2x_1) \sin(3x_2) \exp(-\sqrt{13}x_3)$. 
In the following, $n$ denotes the size of a univariate B-spline basis and $p$ the corresponding polynomial degree, which are common for all spatial dimensions. 

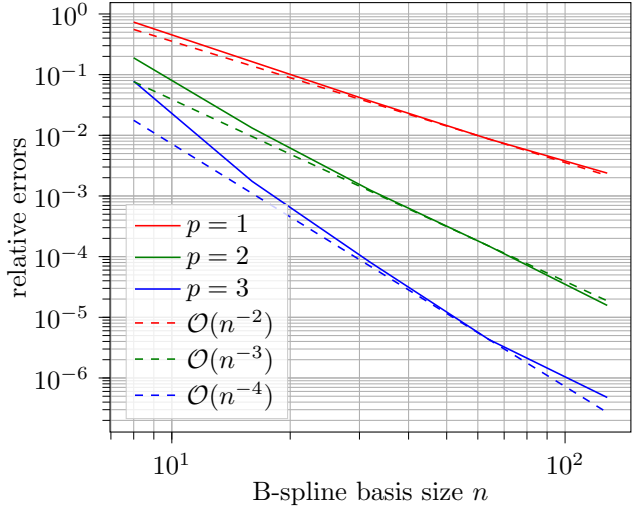
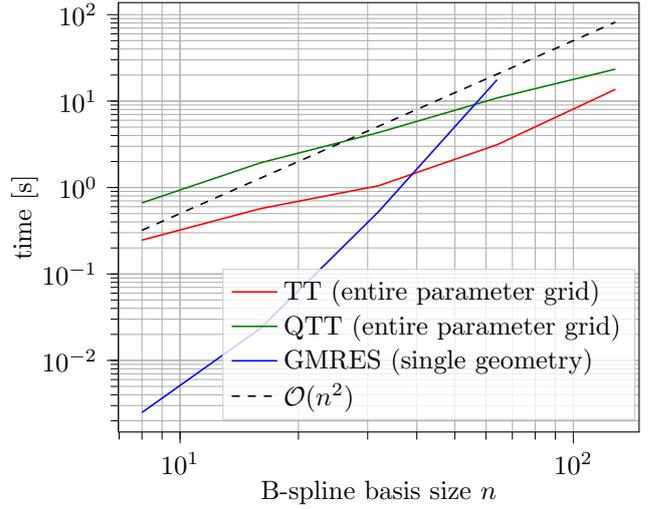
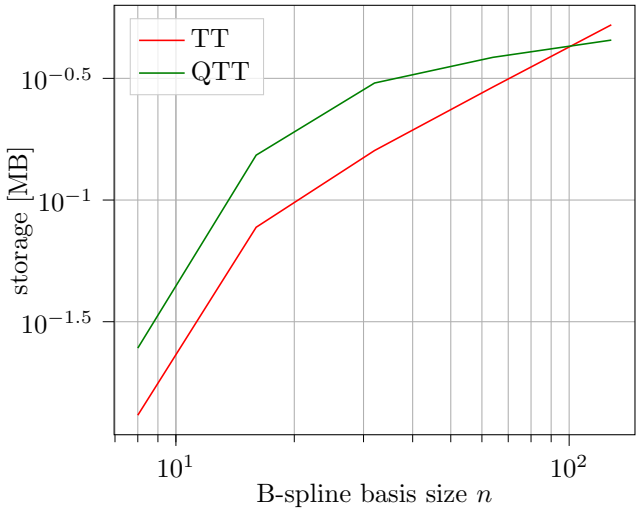
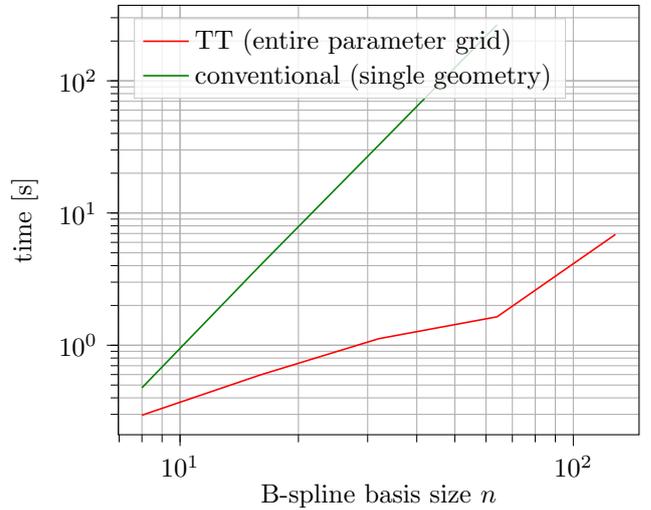
\begin{figure}[t!]
	\centering
	\begin{subfigure}[t]{0.49\textwidth}
\begin{tikzpicture}

\begin{axis}[
legend cell align={left},
legend style={
  fill opacity=0.8,
  draw opacity=1,
  text opacity=1,
  at={(0.03,0.03)},
  anchor=south west,
  draw=white!80!black
},
log basis x={10},
log basis y={10},
tick align=outside,
tick pos=left,
x grid style={white!69.0196078431373!black},
xlabel={B-spline basis size \(\displaystyle n\)},
xmajorgrids,
xmin=6.96440450636899, xmax=147.03338943962,
xminorgrids,
xmode=log,
xtick style={color=black},
y grid style={white!69.0196078431373!black},
ylabel={relative errors},
ymajorgrids,
ymin=1.28454634176838e-07, ymax=1.53912948078408,
yminorgrids,
ymode=log,
ytick style={color=black}
]
\addplot [semithick, red]
table {%
8 0.733709264935134
16 0.16346838603392
32 0.0365903130885607
64 0.00868005787979193
128 0.0023794489886338
};
\addlegendentry{$p=1$}
\addplot [semithick, green!50!black]
table {%
8 0.189093133107667
16 0.0131953411693987
32 0.00124257572488545
64 0.000149450414879074
128 1.57004028629049e-05
};
\addlegendentry{$p=2$}
\addplot [semithick, blue]
table {%
8 0.0780313218290046
16 0.00174558676101176
32 8.02841133315856e-05
64 4.31142576720516e-06
128 4.76628536744321e-07
};
\addlegendentry{$p=3$}
\addplot [semithick, red, dashed]
table {%
8 0.555523704306683
16 0.138880926076671
32 0.0347202315191677
64 0.00868005787979193
128 0.00217001446994798
};
\addlegendentry{$\mathcal{O}(n^{-2})$}
\addplot [semithick, green!50!black, dashed]
table {%
8 0.0765186124180861
16 0.00956482655226076
32 0.00119560331903259
64 0.000149450414879074
128 1.86813018598843e-05
};
\addlegendentry{$\mathcal{O}(n^{-3})$}
\addplot [semithick, blue, dashed]
table {%
8 0.0176595999424723
16 0.00110372499640452
32 6.89828122752825e-05
64 4.31142576720516e-06
128 2.69464110450322e-07
};
\addlegendentry{$\mathcal{O}(n^{-4})$}
\end{axis}

\end{tikzpicture}
		\caption{Solution convergence for degree-$p$ basis, $p\in\left\{1,2,3\right\}$.} \label{fig:conv_err}
	\end{subfigure}
	\hfill
	\begin{subfigure}[t]{0.49\textwidth}
\begin{tikzpicture}

\begin{axis}[
legend cell align={left},
legend style={
  fill opacity=0.8,
  draw opacity=1,
  text opacity=1,
  at={(0.2,0.38)},
  anchor=north west,
  draw=white!80!black
},
log basis x={10},
log basis y={10},
tick align=outside,
tick pos=left,
unbounded coords=jump,
x grid style={white!69.0196078431373!black},
xlabel={B-spline basis size \(\displaystyle n\)},
xmajorgrids,
xmin=6.96440450636899, xmax=147.03338943962,
xminorgrids,
xmode=log,
xtick style={color=black},
y grid style={white!69.0196078431373!black},
ylabel={time [s]},
ymajorgrids,
ymin=0.00148401165875594, ymax=137.783499741107,
yminorgrids,
ymode=log,
ytick style={color=black}
]
\addplot [semithick, red]
table {%
8 0.246262
16 0.569434
32 1.049672
64 3.125921
128 13.682698
};
\addlegendentry{TT (entire parameter grid)}
\addplot [semithick, green!50!black]
table {%
8 0.664193
16 1.928483
32 4.328436
64 10.859625
128 23.425419
};
\addlegendentry{QTT (entire parameter grid)}
\addplot [semithick, blue]
table {%
8 0.002496
16 0.023395
32 0.527651
64 17.640161
128 nan
};
\addlegendentry{GMRES (single geometry)}
\addplot [semithick, black, dashed]
table {%
8 0.32
16 1.28
32 5.12
64 20.48
128 81.92
};
\addlegendentry{$\mathcal{O}(n^{2})$}
\end{axis}

\end{tikzpicture}
		\caption{Solver runtime (quadratic basis). } \label{fig:solver_time}
	\end{subfigure} \\[1em]
	\centering
	\begin{subfigure}[t]{0.49\textwidth}
\begin{tikzpicture}

\begin{axis}[
legend cell align={left},
legend style={
  fill opacity=0.8,
  draw opacity=1,
  text opacity=1,
  at={(0.03,0.97)},
  anchor=north west,
  draw=white!80!black
},
log basis x={10},
log basis y={10},
tick align=outside,
tick pos=left,
x grid style={white!69.0196078431373!black},
xlabel={B-spline basis size \(\displaystyle n\)},
xmajorgrids,
xmin=6.96440450636899, xmax=147.03338943962,
xminorgrids,
xmode=log,
xtick style={color=black},
y grid style={white!69.0196078431373!black},
ylabel={storage [MB]},
ymajorgrids,
ymin=0.0108542784152931, ymax=0.631252353942404,
yminorgrids,
ymode=log,
ytick style={color=black}
]
\addplot [semithick, red]
table {%
8 0.013056
16 0.077312
32 0.159744
64 0.291328
128 0.5248
};
\addlegendentry{TT}
\addplot [semithick, green!50!black]
table {%
8 0.02464
16 0.152928
32 0.302592
64 0.385872
128 0.454176
};
\addlegendentry{QTT }
\end{axis}

\end{tikzpicture}
		\caption{Solution storage needs (quadratic basis).}
		\label{fig:solver_memory}
	\end{subfigure}
	\hfill
	\begin{subfigure}[t]{0.49\textwidth}
\begin{tikzpicture}

\begin{axis}[
legend cell align={left},
legend style={
  fill opacity=0.8,
  draw opacity=1,
  text opacity=1,
  at={(0.03,0.97)},
  anchor=north west,
  draw=white!80!black
},
log basis x={10},
log basis y={10},
tick align=outside,
tick pos=left,
unbounded coords=jump,
x grid style={white!69.0196078431373!black},
xlabel={B-spline basis size \(\displaystyle n\)},
xmajorgrids,
xmin=6.96440450636899, xmax=147.03338943962,
xminorgrids,
xmode=log,
xtick style={color=black},
y grid style={white!69.0196078431373!black},
ylabel={time [s]},
ymajorgrids,
ymin=0.210445817165706, ymax=372.507358532735,
yminorgrids,
ymode=log,
ytick style={color=black}
]
\addplot [semithick, red]
table {%
8 0.295649
16 0.596063
32 1.120147
64 1.640698
128 6.892718
};
\addlegendentry{TT (entire parameter grid)}
\addplot [semithick, green!50!black]
table {%
8 0.476336
16 4.019093
32 32.580096
64 265.15434
128 nan
};
\addlegendentry{conventional (single geometry)}
\end{axis}

\end{tikzpicture}
		\caption{Stiffness assembly time (quadratic basis).}
		\label{fig:assembly}
	\end{subfigure}
	\caption{Test case 1: Solution convergence, solver runtime, memory requirements for storing the solution, and assembly time using IGA solvers based on TT, QTT, and GMRES. The TT and QTT-based IGA solvers resolve the parameter-dependent cylinder's geometry, while the GMRES is used only for a single geometry configuration (corresponding to the largest element from the parameter grid). The parameter-dependence is resolved using a collocation grid of constant size $\ell=8$.}
	\label{fig:convergence2}
\end{figure}

As a first study, the parameter dependence is resolved using a collocation grid of constant size $\ell=8$.
The corresponding results are presented in Figure \ref{fig:convergence2}. 
In Figure~\ref{fig:conv_err}, the convergence of the TT solver is shown for gradually refined linear, quadratic, and cubic B-spline bases.
As expected from theory, the convergence order is $\mathcal{O}(n^{-(p+1)})$.  
For the remaining results, a quadratic B-spline basis is used, however, we note that the results for linear and cubic bases are very similar to the presented ones.
In Figure~\ref{fig:solver_time}, the TT and QTT formats are compared against one another in terms of computation time, where it is found that the TT solver has a slight advantage. 
In both cases, the TT-rank stagnates after a certain basis size. 
The complexity of both formats increases in similar fashion as the B-spline basis is refined.
The scaling of a standard GMRES solver is additionally shown. In this case, the parameter-dependent deformation is omitted and only a single parameter value is considered. 
Nevertheless, it can be observed that the runtime of the GMRES solver increases much faster than that of the TT/QTT-IGA solvers, even though in the latter case the parameter dependence is taken into account as well. 
For a fairer comparison, the runtime of the GMRES solver must be multiplied by a factor of $\ell=8$, which is the size of the parameter grid employed in the TT/QTT solvers.
Regarding implementation specifics, it should be noted that neither the TT/QTT nor the GMRES solvers use preconditioning, however, the local systems solved during the AMEn iterations benefit from Jacobi preconditioning.
In Figure~\ref{fig:solver_memory}, the TT and QTT formats are compared with respect to memory requirements.
It can be observed that the QTT format becomes the most memory-efficient option after $n=60$ basis functions per spatial dimension. 
Finally, Figure~\ref{fig:assembly} shows how the TT-based assembly significantly outperforms the conventional IGA assembly.  
Note that the assembly of the parameter-dependent TT-operators is orders of magnitude faster, even though the conventional approach is applied for a single parameter value only. 
The complexity of the assembly is reduced to $\mathcal{O}(n^2)$ from the complexity of $\mathcal{O}(n^3)$ corresponding to the conventional IGA stiffness assembler.

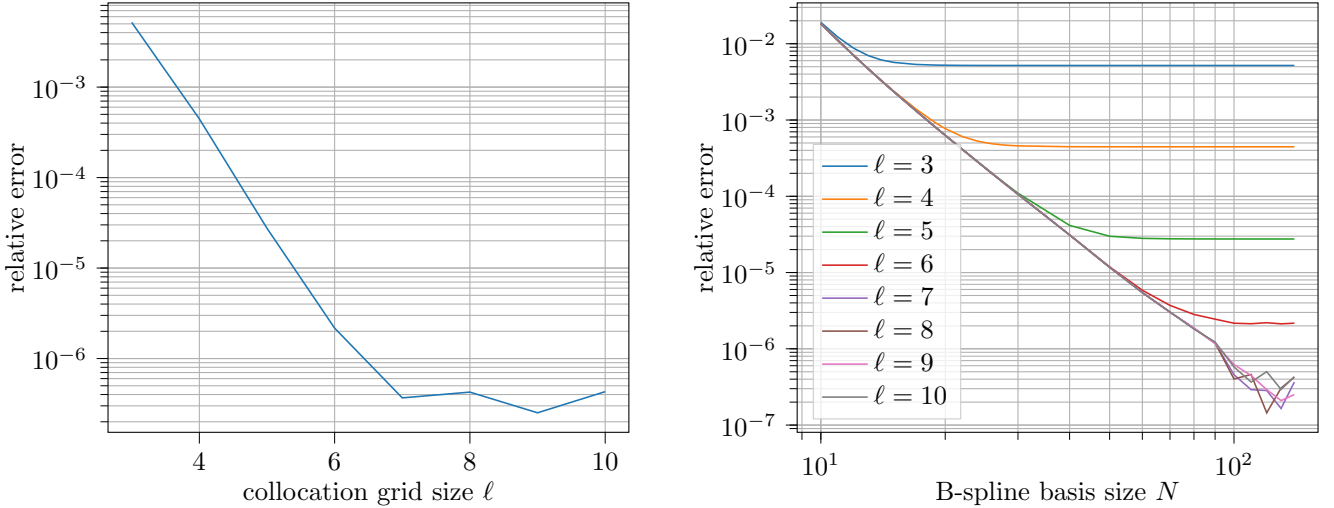
\begin{figure}[t!]
	\centering
	\begin{subfigure}[t]{0.49\textwidth}
\begin{tikzpicture}

\definecolor{color0}{rgb}{0.12156862745098,0.466666666666667,0.705882352941177}

\begin{axis}[
log basis y={10},
tick align=outside,
tick pos=left,
x grid style={white!69.0196078431373!black},
xlabel={collocation grid size \(\displaystyle \ell\)},
xmajorgrids,
xmin=2.65, xmax=10.35,
xminorgrids,
xtick style={color=black},
y grid style={white!69.0196078431373!black},
ylabel={relative error},
ymajorgrids,
ymin=1.52988271752459e-07, ymax=0.00852669760784667,
yminorgrids,
ymode=log,
ytick style={color=black}
]
\addplot [semithick, color0]
table {%
3 0.00518857394342701
4 0.000445497302972593
5 2.75341094714259e-05
6 2.16845748255872e-06
7 3.67465063405614e-07
8 4.25404316454838e-07
9 2.51414886827013e-07
10 4.30459170988094e-07
};
\end{axis}

\end{tikzpicture}
		\caption{Solution convergence with respect to the size of the collocation grid $\ell$ for a cubic B-spline basis of fixed size ($N=128$).}\label{fig:ell1}
	\end{subfigure}%
	\hfill
	\begin{subfigure}[t]{0.49\textwidth}
\begin{tikzpicture}

\definecolor{color0}{rgb}{0.12156862745098,0.466666666666667,0.705882352941177}
\definecolor{color1}{rgb}{1,0.498039215686275,0.0549019607843137}
\definecolor{color2}{rgb}{0.172549019607843,0.627450980392157,0.172549019607843}
\definecolor{color3}{rgb}{0.83921568627451,0.152941176470588,0.156862745098039}
\definecolor{color4}{rgb}{0.580392156862745,0.403921568627451,0.741176470588235}
\definecolor{color5}{rgb}{0.549019607843137,0.337254901960784,0.294117647058824}
\definecolor{color6}{rgb}{0.890196078431372,0.466666666666667,0.76078431372549}

\begin{axis}[
legend cell align={left},
legend style={
  fill opacity=0.8,
  draw opacity=1,
  text opacity=1,
  at={(0.03,0.03)},
  anchor=south west,
  draw=white!80!black
},
log basis x={10},
log basis y={10},
tick align=outside,
tick pos=left,
x grid style={white!69.0196078431373!black},
xlabel={B-spline basis size \(\displaystyle N\)},
xmajorgrids,
xmin=8.76382301087964, xmax=159.747635051735,
xminorgrids,
xmode=log,
xtick style={color=black},
y grid style={white!69.0196078431373!black},
ylabel={relative error},
ymajorgrids,
ymin=8.0381002668518e-08, ymax=0.0343494962647715,
yminorgrids,
ymode=log,
ytick style={color=black}
]
\addplot [semithick, color0]
table {%
10 0.0190535726427983
11 0.012200675136364
12 0.00875227987984051
13 0.00701423755648422
14 0.00614129837509262
15 0.00570054934041328
17 0.00535264246432406
19 0.00524842533706019
20 0.00522632612892186
22 0.00520462918213292
24 0.00519598335229361
26 0.00519219722543574
28 0.00519045436369137
30 0.00518958857855809
40 0.00518864592779735
50 0.00518857880052288
60 0.00518857317464226
70 0.00518857327085962
80 0.00518857267980836
90 0.00518857327895434
100 0.00518857387790987
110 0.00518857374870579
120 0.0051885739069046
130 0.00518857407440378
140 0.00518857394342701
};
\addlegendentry{$\ell=3$}
\addplot [semithick, color1]
table {%
10 0.01834065176988
11 0.0110536674687974
12 0.00706435380549308
13 0.00474286334967951
14 0.00331769472300094
15 0.00240521855378163
17 0.00139123365925371
19 0.000909989432439946
20 0.000771804146095667
22 0.000606677164640404
24 0.000526586938590757
26 0.000487269717915131
28 0.000467937853941642
30 0.000457974150256403
40 0.0004465754459339
50 0.000445648547867143
60 0.000445529561732103
70 0.000445506479713541
80 0.000445500756030056
90 0.000445498755180448
100 0.000445497458675533
110 0.000445497245284279
120 0.000445497465987284
130 0.000445497308284692
140 0.000445497302972593
};
\addlegendentry{$\ell=4$}
\addplot [semithick, color2]
table {%
10 0.0183355252959807
11 0.0110450412719462
12 0.0070505707464333
13 0.00472216637960654
14 0.00328793842888912
15 0.00236392209717774
17 0.00131840731014138
19 0.00079408996267959
20 0.000630976689777935
22 0.000412851405153689
24 0.00028222144108079
26 0.000199415079185392
28 0.000145900272021426
30 0.000109786627767592
40 4.15440030863898e-05
50 2.99148023638423e-05
60 2.80686431284911e-05
70 2.77009866060453e-05
80 2.75912157526081e-05
90 2.75607471622095e-05
100 2.75416065004638e-05
110 2.75335523326232e-05
120 2.75381717956299e-05
130 2.75352202600464e-05
140 2.75341094714259e-05
};
\addlegendentry{$\ell=5$}
\addplot [semithick, color3]
table {%
10 0.0183355048380147
11 0.011045008334967
12 0.00705051853556251
13 0.00472208747689521
14 0.00328782460489024
15 0.00236376342867809
17 0.00131812211286759
19 0.000793616039581827
20 0.000630379760659895
22 0.000411938502065254
24 0.00028088407020133
26 0.000197517137167644
28 0.000143295426020559
30 0.000106299225961652
40 3.11822322478536e-05
50 1.18938064507272e-05
60 5.86076913809266e-06
70 3.70983815178227e-06
80 2.82200581176825e-06
90 2.44807462133312e-06
100 2.16497445901578e-06
110 2.13163127677624e-06
120 2.20003790428537e-06
130 2.12413111060518e-06
140 2.16845748255872e-06
};
\addlegendentry{$\ell=6$}
\addplot [semithick, color4]
table {%
10 0.018335504734645
11 0.0110450081087852
12 0.00705051820436482
13 0.00472208700234253
14 0.0032878239321886
15 0.0023637624707199
17 0.00131812016205959
19 0.000793613069637817
20 0.000630375974881712
22 0.0004119330179095
24 0.000280875650545843
26 0.00019750527483672
28 0.000143278873322375
30 0.000106277269531938
40 3.1106931899098e-05
50 1.17011556250124e-05
60 5.4717884775978e-06
70 3.02948529637141e-06
80 1.86283972174626e-06
90 1.18084120488499e-06
100 4.57217416795431e-07
110 2.91026750340172e-07
120 2.84021241878804e-07
130 1.65097874526101e-07
140 3.67465063405614e-07
};
\addlegendentry{$\ell=7$}
\addplot [semithick, color5]
table {%
10 0.0183355047274388
11 0.0110450081104222
12 0.0070505182115901
13 0.0047220870055129
14 0.00328782392216436
15 0.0023637624563904
17 0.00131812027999122
19 0.000793613194399898
20 0.000630375926335698
22 0.00041193273943058
24 0.000280875562587432
26 0.000197504515820648
28 0.000143279585564476
30 0.000106277095933416
40 3.11076865359043e-05
50 1.17021035219725e-05
60 5.45527015655512e-06
70 3.02806102658358e-06
80 1.8493028095846e-06
90 1.19841185728245e-06
100 4.02409330825723e-07
110 4.605005596793e-07
120 1.44909671413484e-07
130 3.05770773677636e-07
140 4.25404316454838e-07
};
\addlegendentry{$\ell=8$}
\addplot [semithick, color6]
table {%
10 0.0183355047336093
11 0.011045008117743
12 0.00705051823445606
13 0.00472208698491367
14 0.00328782389670345
15 0.00236376241590068
17 0.00131812035334204
19 0.000793613206621881
20 0.000630375975014099
22 0.000411932887429672
24 0.000280875645568915
26 0.000197505484047787
28 0.000143279693079726
30 0.000106276867227193
40 3.11069016609935e-05
50 1.1716799321394e-05
60 5.47964459902139e-06
70 3.01745279581213e-06
80 1.87704370105379e-06
90 1.16854771318364e-06
100 6.20905709437062e-07
110 4.46077249255173e-07
120 2.91598715890771e-07
130 2.09269819018235e-07
140 2.51414886827013e-07
};
\addlegendentry{$\ell=9$}
\addplot [semithick, white!49.8039215686275!black]
table {%
10 0.018335504790666
11 0.0110450081842281
12 0.00705051828099449
13 0.00472208707119395
14 0.00328782403679395
15 0.00236376251902509
17 0.00131812030068397
19 0.00079361325282993
20 0.000630376102554388
22 0.000411932837385335
24 0.00028087586798919
26 0.000197505196860213
28 0.000143279433006003
30 0.000106279163550296
40 3.11066326331457e-05
50 1.17139112442372e-05
60 5.45906045542481e-06
70 3.03793621176855e-06
80 1.81347429220968e-06
90 1.22729350735728e-06
100 5.78262791925652e-07
110 3.66645486034488e-07
120 5.01557966208225e-07
130 2.93519349350809e-07
140 4.30459170988094e-07
};
\addlegendentry{$\ell=10$}
\end{axis}

\end{tikzpicture}
		\caption{Solution convergence with respect to the size of the cubic B-spline basis $N$ and the size of the collocation grid $\ell$.}\label{fig:ell2}
	\end{subfigure}
	\caption{Test case 1: Solution convergence for collocation grids of increasing size.}
	\label{fig:ell}
\end{figure}

Using the same problem setting, a second study is performed to examine the convergence of the parameter-dependent solution with respect to the number of collocation points $\ell$. 
The results are presented in Figure \ref{fig:ell}, where an exponential decrease of the error is observed for an increasing size of the collocation grid.
In this case, a cubic B-spline basis of fixed size $n=128$ is used.
It is also evident that the discretization of the physical domain becomes the limiting factor in terms of solution accuracy for $\ell \geq 7$.  
The same effect is further illustrated in Figure \ref{fig:ell2}, where the size of the cubic B-spline basis is increased along with the size of the collocation grid.

\subsection{Test case 2: TT-IGA solver performance for multiple parameter dependencies}
\label{subsec:test-case-2}
\begin{figure}[t!]
	\centering
	\begin{subfigure}[t]{0.49\textwidth}
		\centering
		\includegraphics[width=\textwidth]{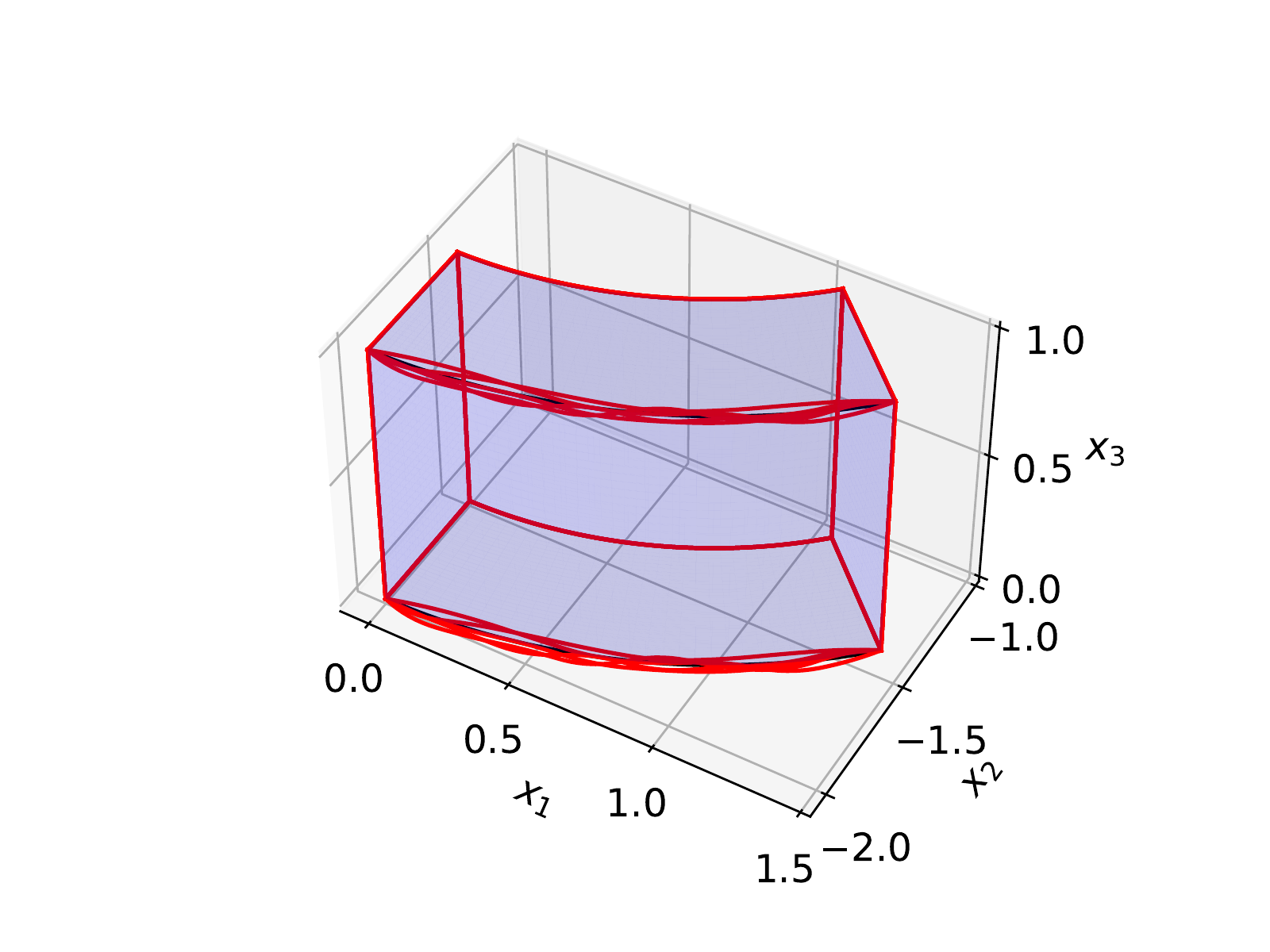}
		\caption{The C-shaped domain for different parameter combinations ($N_{\text{p}}=8$). } \label{fig:C8_geom}
	\end{subfigure}
	\hfill
	\begin{subfigure}[t]{0.49\textwidth}
		\centering
		\includegraphics[width=\textwidth]{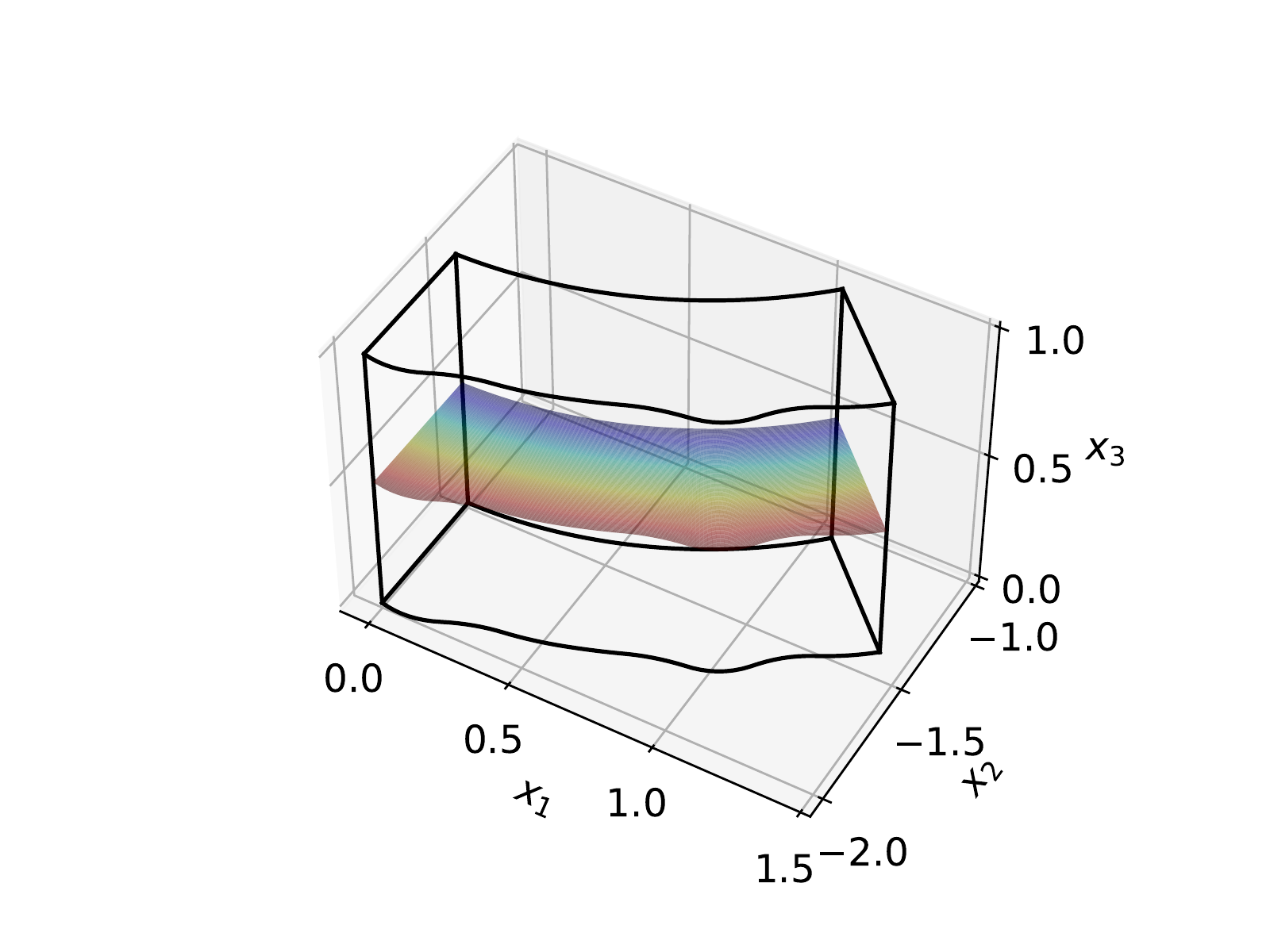}
		\caption{Solution in the $x_1 x_2$-plane for $x_3=0.5$. A zero Dirichlet condition is enforced along the plane defined by the inner and outer radius and a zero Neumann boundary condition in the remaining boundaries.} \label{fig:C8_sol}
	\end{subfigure}
	\caption{Test case 2: Quarter of a C-shaped domain with parameter-dependent deformations.}
	\label{fig:C8}
\end{figure}

We now assess the performance of the TT-IGA solver for an increasing number of parameters affecting the geometry of the physical domain. 
To that end, a quarter of a C-shaped domain with a perturbed outer radius is considered, as shown in Figure~\ref{fig:C8_geom}. 
The inner radius is $r_{\text{in}}=1.5$, while the outer radius is given by
\begin{align}
	r_{\text{out}}(\alpha,\btheta) = 2+ \sum\limits_{k=1}^{N_\text{p}} \theta_k b^{(\text{out})}_k(\alpha\frac{\pi}{4}),
\end{align}
where $\alpha\in(0,\pi/4)$ is the angle spanning the domain, $N_\text{p}$ is the number of parameters, $\left\{b^{(\text{out})}_k\right\}_{k=1}^{N_\text{p}}$ is a B-spline basis of size $N_\text{p}$, and $\theta_k \in [-0.0.5,0.05]$. 
The parametrization of the domain $\Omega(\btheta)$ is given by
\begin{align}
	G(\bm{y},\btheta) = \begin{pmatrix}
		(1.5+y_1 r_{\text{out}}(y_2,\btheta))\cos(y_2\pi/4)  \\
		(1.5+y_1 r_{\text{out}}(y_2,\btheta))\sin(y_2\pi/4) \\
		y_3
	\end{pmatrix},\quad \bm{y}\in[0,1]^3.
	\end{align}
The following equation with Dirichlet conditions on the exterior and interior circular surfaces and Neumann conditions on the remaining boundaries is solved on the parameter-dependent geometry described above:
\begin{align}
	\Delta u(\cdot,\btheta)=0&,\quad \text{in } \Omega(\btheta),\\
	\partial_{\bm{\nu}}u = 0& ,\quad \text{on } \partial\Omega(\btheta)\cap\{G(\bm{y},\btheta) \: : \: y_1\in(0,1)\}, \\
	u = 1& ,\quad \text{on } \partial\Omega(\btheta)\cap\{G(\bm{y},\btheta) \: : \: y_1=1\}, \\
	u = 0& ,\quad \text{on } \partial\Omega(\btheta)\cap\{G(\bm{y},\btheta) \: : \: y_1=0\}.
\end{align}
See Figure \ref{fig:C8_sol} for a representation of the solution for a specific parameter combination.

In the following, the spatial discretization is based on quadratic B-splines with a basis size $\bm{n} = (40,20,80)$, while the parameter grid is constructed with $\ell_1 = \cdots =  \ell_{N_{\text{p}}} = \ell = 8$ nodes per parameter. 
For the construction of the stiffness matrix, the QTT format is used in order to limit the memory consumption of the pointwise multiplications in the TT-format. 
The number of parameters is varied and the results are reported in Table \ref{tab:iga_multiparam}. 
Regarding time and computational complexity, a slight exponential scaling with respect to the number of parameters is observed $(\sim 1.9^{N_\text{p}})$, therefore, the TT-IGA solver remains applicable up to a moderate number of parameters.
With respect to storage needs, the memory requirements for the stiffness operator and the system matrix are two orders of magnitude higher than for the solution tensor. 
This can be further reduced if the band diagonal structure of the cores is exploited by storing them in sparse format. 
The error values have been verified for several geometry configurations in order to make sure that the error does not grow while increasing the number of parameters.

\begin{table}[t!]
	\centering
	\begin{tabular}{   c c   c  c  c    } \hline \hline
		$N_{\text p}$ & Stiffness assembly [s]  & Solver runtime [s] & Operator storage [MB]    & Solution storage [MB]   \\
		\hline
		2 & $15.4$ & $2.2$    &$14$ & $0.08$ \\ 
		3 & $28.2$ & $6.2$   &$31$ & $0.14$ \\ 
		4 & $46.4$ & $15.4$   &$49$ & $0.23$ \\ 
		5 & $68.7$ & $38.1$  &$72$ & $0.39$ \\ 
		6 & $92$   & $87.3$  &$94$ & $0.6$ \\ 
		7 & $148$  & $190$     &$134$ & $0.88$ \\ 
		8 & $238$  & $355$    &$165$ & $1.3$ \\ 
		9 & $331$  & $598$    &$202$ & $1.8$ \\ 
		10 & $426$ & $1012$   &$245$ & $2.4$  \\ \hline \hline
	\end{tabular}
	\caption{Test case 2: Computational and storage costs for assembling the stiffness matrix and solving the system in the TT-format for an increasing number of parameters $N_{\text p}$.}
	\label{tab:iga_multiparam}
\end{table}

\subsection{Test case 3: Domain with material discontinuity}
In the next numerical investigation we consider the case where the coefficient function $\kappa$ in the model BVP \eqref{eq:bvp} is piecewise discontinuous with respect to the coordinates of the physical domain. 
Exemplarily, this case becomes of great interest when modeling domains comprising materials with different properties. 
We restrict ourselves to the case where the material coefficient represented in the reference domain $\hat{\kappa}$ is smooth over Cartesian partitions of the reference domain $[0,1]^3$. 
The discontinuity of $\hat{\kappa}$ along the subdomain boundaries implies lack of smoothness for the solution along those surfaces. 
The B-spline basis needs to be chosen accordingly by inserting the bounds of the Cartesian partition in the knots vector with increased multiplicity, that is, the points need to appear $p$ times, where $p$ is the degree of the B-splines, also see Section~\ref{subsec:cad}. 

\begin{figure}[t!]
	\centering
	\begin{subfigure}[t]{0.49\textwidth}
		\centering
		\includegraphics[width=\textwidth]{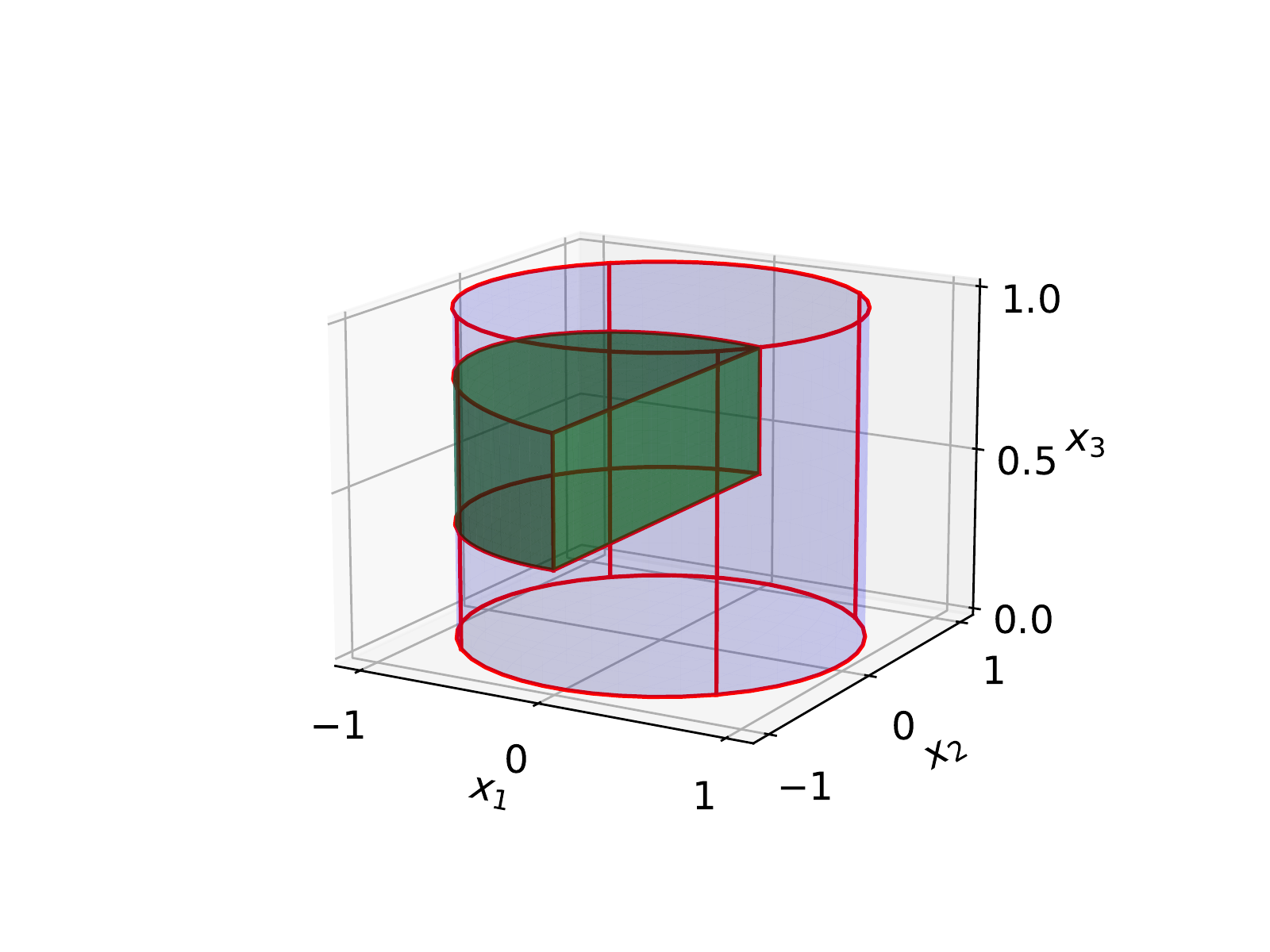}
		\caption{Computational domain for $\theta_2=-0.05,\:\theta_3 = \theta_4 = 0.05$.}
	\end{subfigure}
	\hfill
	\begin{subfigure}[t]{0.49\textwidth}
		\centering
		\includegraphics[width=\textwidth]{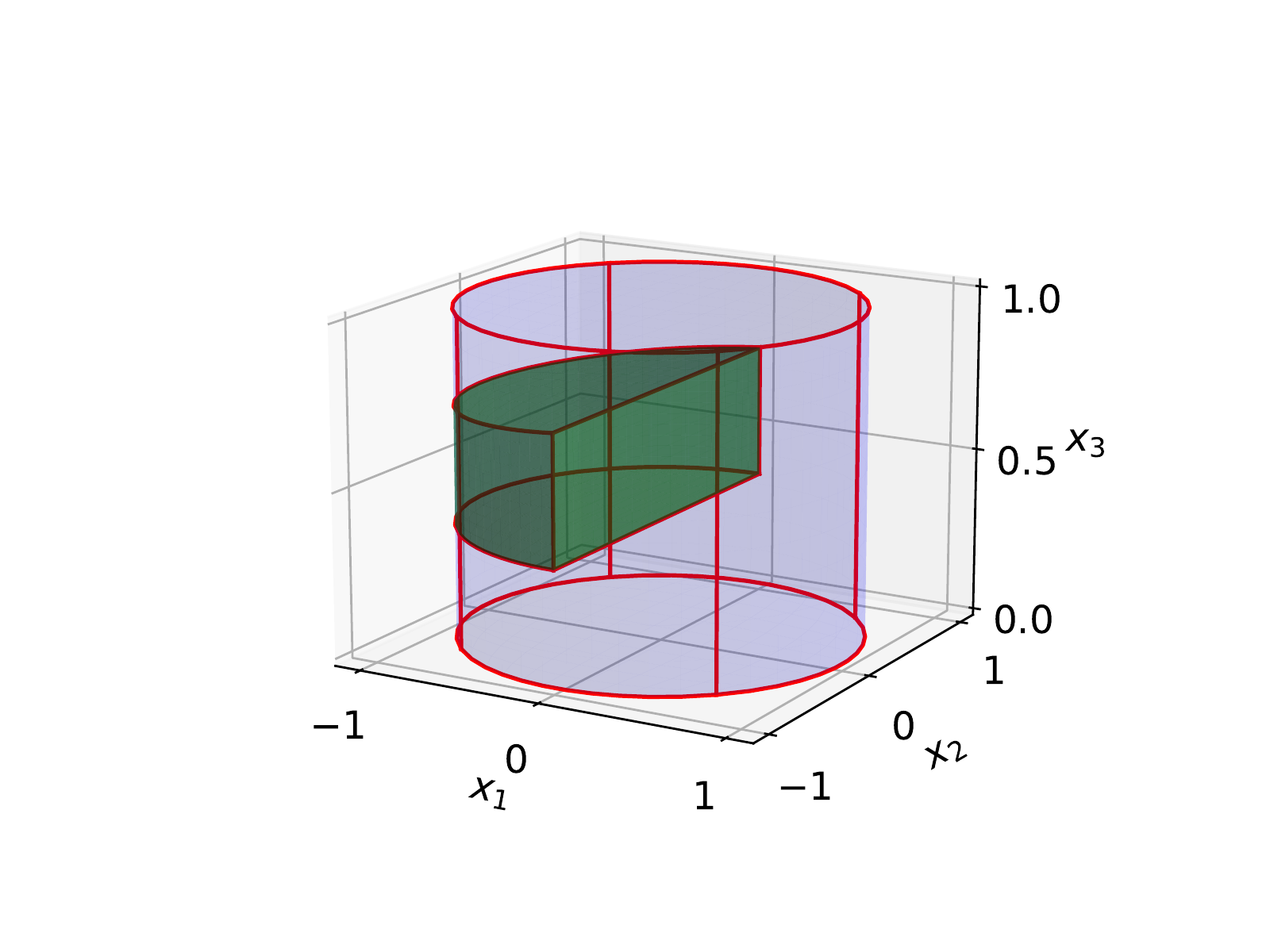}
		\caption{Computational domain for $\theta_2=0.05,\:\theta_3 = \theta_4 = 0.05$.}
	\end{subfigure}
	\caption{Test case 3: Cylindrical domain with material jump. The material coefficient is $\kappa=5+0.5\theta_3,\: \theta_3\in[-1,1]$, inside the subdomain corresponding to the material filling (denoted with green) and $\kappa=1$ otherwise.}
	\label{fig:cylinder_material}
\end{figure}

In particular, we consider a cylindrical domain with material filling, as shown in Figure~\ref{fig:cylinder_material}. 
The cylinder has radius $r=1$ and height $h=1$ along the $x_3$-coordinate, such that $\Omega = \{(x_1,x_2,x_3)\: : \: x_1^2+x_2^2<1 \: \wedge \: 0<x_3<1 \}$. 
The parameter-dependent BVP to be solved is the Laplace equation with Dirichlet boundary conditions applied on the top $(x_3=1)$ and the bottom $(x_3=0)$ boundaries of the cylinder, while Neumann boundary conditions are applied on the remaining boundary. 
The BVP reads
\begin{subequations}
\label{eq:laplace-discontinuity}
\begin{align}
	\nabla\cdot(\kappa(\cdot,\btheta) \nabla u(\cdot,\btheta)) = 0,& \quad \text{in }\Omega , \\
	u(\cdot,\btheta)) = 10,&\quad \text{on } \partial\Omega \cap \{(x_1,x_2,x_3)\: : \: x_3 = 0\}, \\
	u(\cdot,\btheta)) = 0,&\quad \text{on } \partial\Omega \cap \{(x_1,x_2,x_3)\: : \: x_3 = 1\}, \\
	\partial_{\bm{\nu}} u(\cdot,\btheta)) = 0,&\quad \text{on } \partial\Omega \cap \{(x_1,x_2,x_3)\: : \: 0 < x_3 < 1\}.
\end{align}
\end{subequations} 
The discontinuous material coefficient is given as
\begin{align}
	\kappa(\bm{x},\btheta) = \begin{cases}
		5+0.5\theta_1,&\quad \bm{x} \in \Upsilon(\theta_2,\theta_3,\theta_4), \\ 
		1 ,&  \bm{x} \in \Omega \setminus \Upsilon(\theta_2,\theta_3,\theta_4),
	\end{cases}
\end{align} 
where the parameter-dependent subdomain is defined as
\begin{equation} 
\Upsilon(\theta_2,\theta_3,\theta_4)  =\{ (x_1,x_2,x_3) \: : \: x_1^2+x_2^2<1 \: \wedge \: x_1<0 \: \wedge \: 0.3+\theta_4 < x_3 < 0.7+\theta_3 + \theta_2 x_1 \},  
\end{equation}	
and depends on 3 parameters, see Figure \ref{fig:cylinder_material}.
The parameter support is chosen as $\Xi = [-1,1] \times [-0.05,0.05] \times [-0.05,0.05] \times [-0.05,0.05]$.

\begin{figure}[t!]
	\centering
	\begin{subfigure}[b]{0.49\textwidth}
		\centering
		\includegraphics[width=\textwidth]{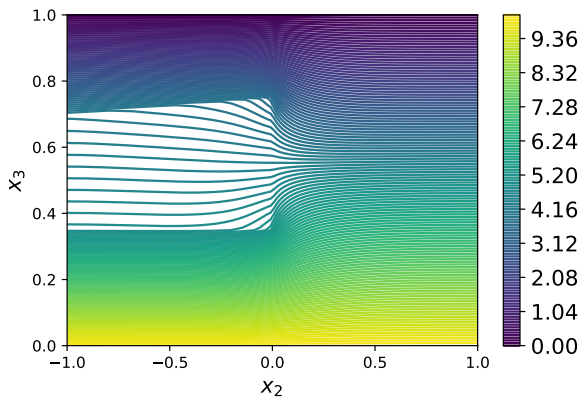}
		\caption{Solution of the TT-IGA solver.}
	\end{subfigure}
	\hfill
	\begin{subfigure}[b]{0.49\textwidth}
		\centering
		\includegraphics[width=\textwidth]{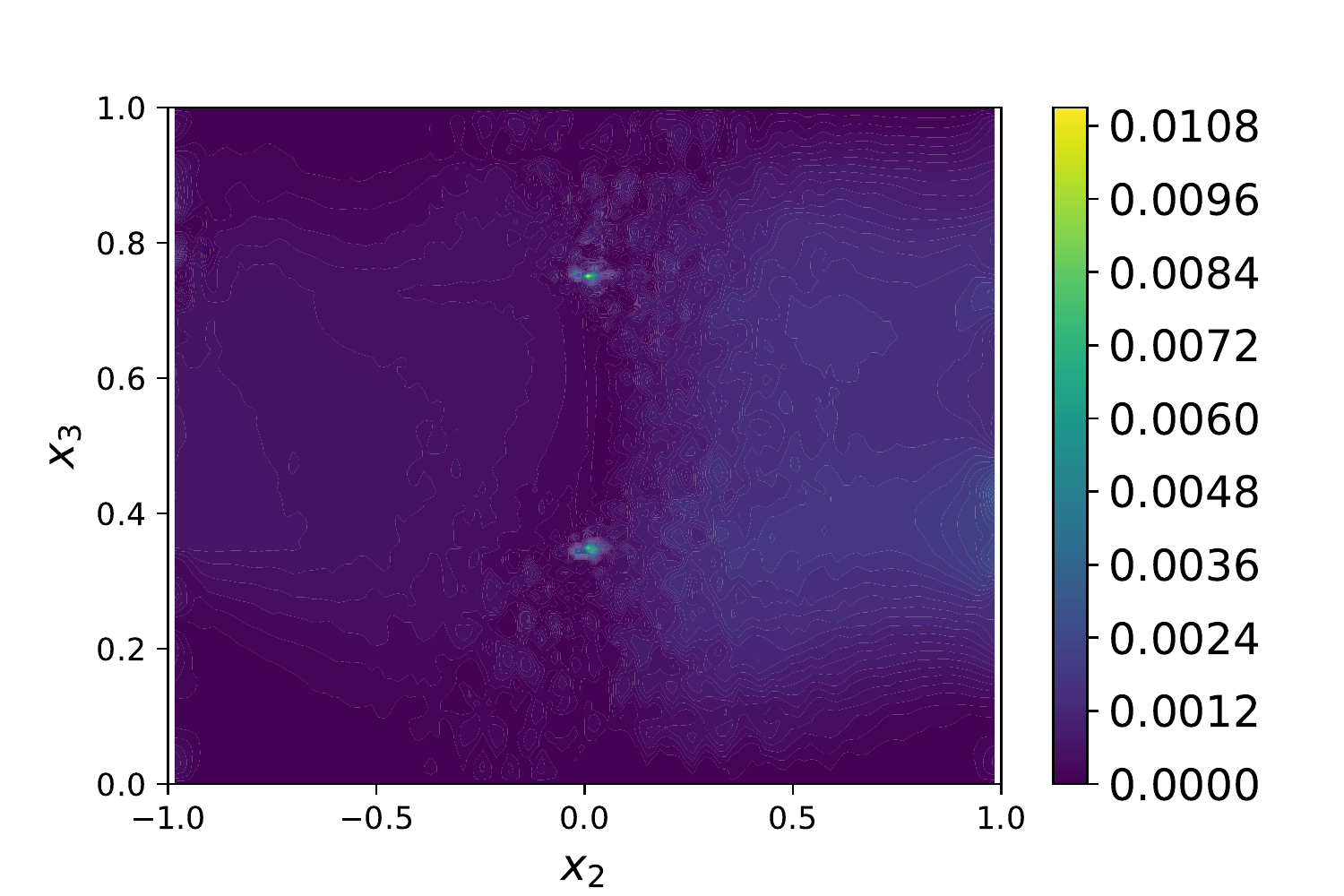}
		\caption{Pointwise absolute error compared to FEM solution.}
	\end{subfigure}
	\centering
	\caption{Test case 3: Solution to the problem \eqref{eq:laplace-discontinuity} along the plane $x_2=0$ for $\theta_1 = \theta_2 = \theta_3 = \theta_4=0.05$ and pointwise absolute error with respect to an FEM reference solution.}
	\label{fig:jump_solution}
\end{figure}

In Figure~\ref{fig:jump_solution}, the solution to problem \eqref{eq:laplace-discontinuity} for $\theta_1 = \theta_2 = \theta_3 = \theta_4 = 0.05$ is plotted along the plane $x_2=0$. 
The size of the B-spline basis is $\bm{n}=(80,80,80)$, while that of the collocation grid is $\bm{\ell}=(12,12,12,12)$. 
The relative accuracy of the AMEn solver is set to $\epsilon=10^{-6}$. 
The IGA solution is compared against a high fidelity FEM solution obtained using the \texttt{FEniCS} software \cite{alnaes2015fenics}. 
The FEM solver needs approximately $41$~s for a single solution, equivalently, for a single realization of the parameter vector. 
In comparison, the runtime of TT-IGA solver is approximately equal to $23$~s for an approximation over the entire parameter space. 
Accordingly, the number of floating point entries needed to store a single FEM solution is approximately equal to $293000$, while the TT-based solution can be stored using merely $107000$ entries for the combined state-parameter approximation. 
While discretizing the physical space, i.e. for increasing values $\bm{n}$, the TT rank of the solution tensor as well as the TT rank of the discrete operator increase up to a plateau, as can be seen in Table~\ref{tab:iga_gmres}).




Furthermore, the TT-IGA solver is compared against a classical GMRES solver in terms of computation time, see Table~\ref{tab:iga_gmres}.
Note that the TT-IGA solver addresses the parameter-dependent problem, while the GMRES solver is employed for a fixed geometry. 
In the former case, the parameter dependence is resolved using a collocation grid of size $\bm{\ell} = \left(12,12,12,12\right)$. 
For coarse spatial discretizations, the runtime of the GMRES solver is fast enough to allow the consideration of other techniques to resolve the parametric dependency, e.g. using stochastic collocation on sparse grids. 
However, for finer discretizations, equivalently, for higher $\bm{n}$ values, solving the TT system for the entire tensor-product grid is faster even if compared to solving the system with the GMRES solver for a single parameter realization, equivalently, geometry configuration. 

\begin{table}[b!]
	\centering
	\begin{tabular}{   c   c  c  c  c  c  } \hline \hline
		B-spline basis terms per dimension & 20 & 30 & 40 & 50 & 60  \\ \hline
		TT-IGA runtime [s] & $3.6$ & $6.7$ & $13.2$ & $22.6$ & $33.1$ \\ 
		GMRES runtime [s] & $0.21$ & $3.4$ & $21.7$ & $85$ & $298$  \\ 
			mean TT rank, solution tensor & $28.625$ & $29.875$ & $33.125$ & $32.375$ & $32.375$  \\ 
		mean TT rank, tensor-operator & $20.750$ & $21.625$ & $22.125$ & $22.250$ & $25.500$ \\ \hline\hline
	\end{tabular}
	\caption{Test case 3: Runtime and rank comparison between the TT-IGA AMEn solver for the entire parameter grid and the GMRES solver for a fixed geometry, for an increasingly refined spatial discretization.}
	\label{tab:iga_gmres}
\end{table}

\subsection{Test case 4: Waveguide problem}
\begin{figure}[t!]
	\centering
	\begin{subfigure}[t]{0.49\textwidth}
		\centering
		\includegraphics[width=\textwidth]{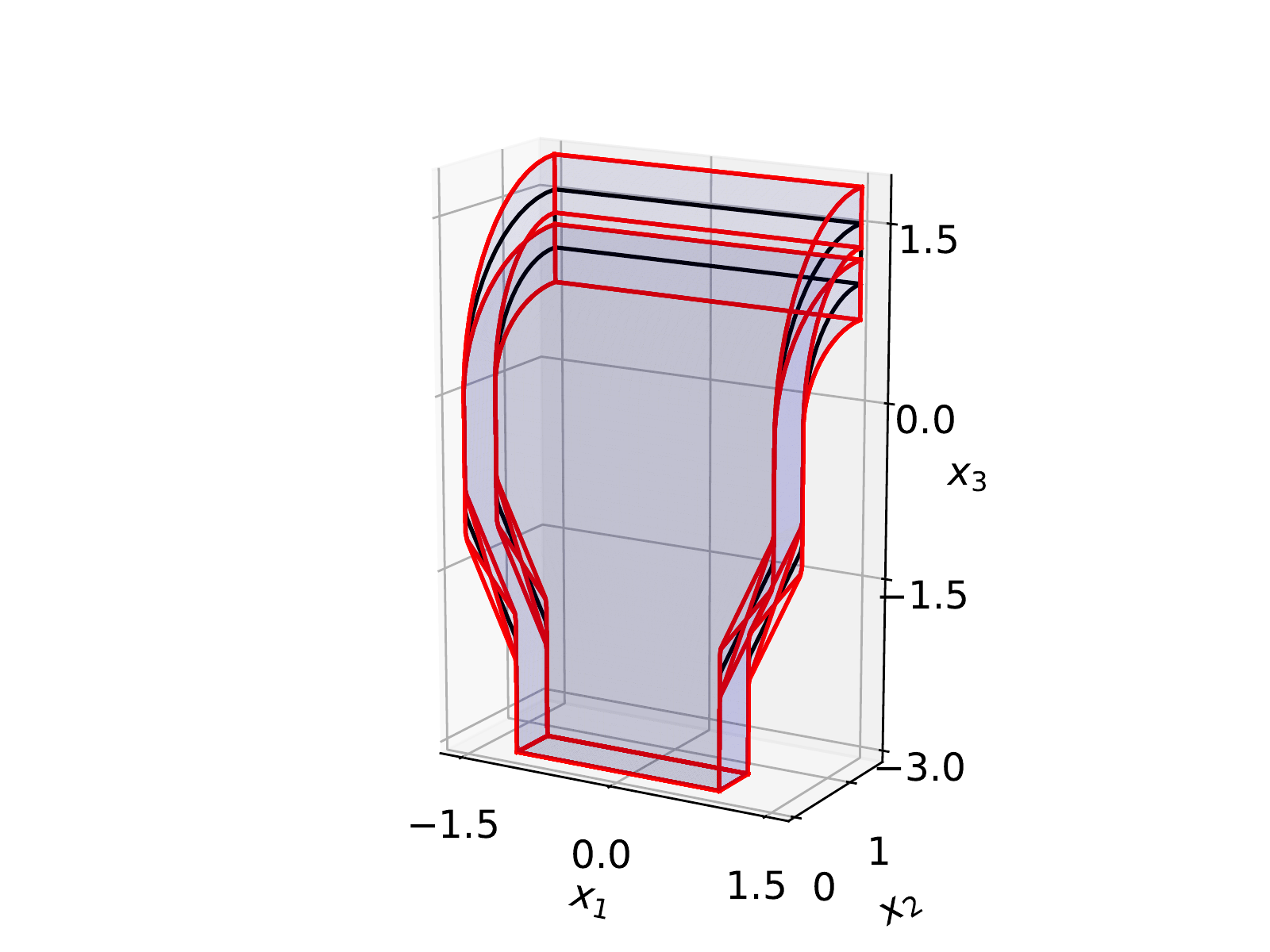}
		\caption{Waveguide geometry configurations. The black contour corresponds to the nominal configuration. Geometry configurations that correspond to the vertices of the parameters' support space are shown in red.} \label{fig:wg_parms}
	\end{subfigure}
	\hfill
	\begin{subfigure}[t]{0.49\textwidth}
		\centering
		\includegraphics[width=\textwidth]{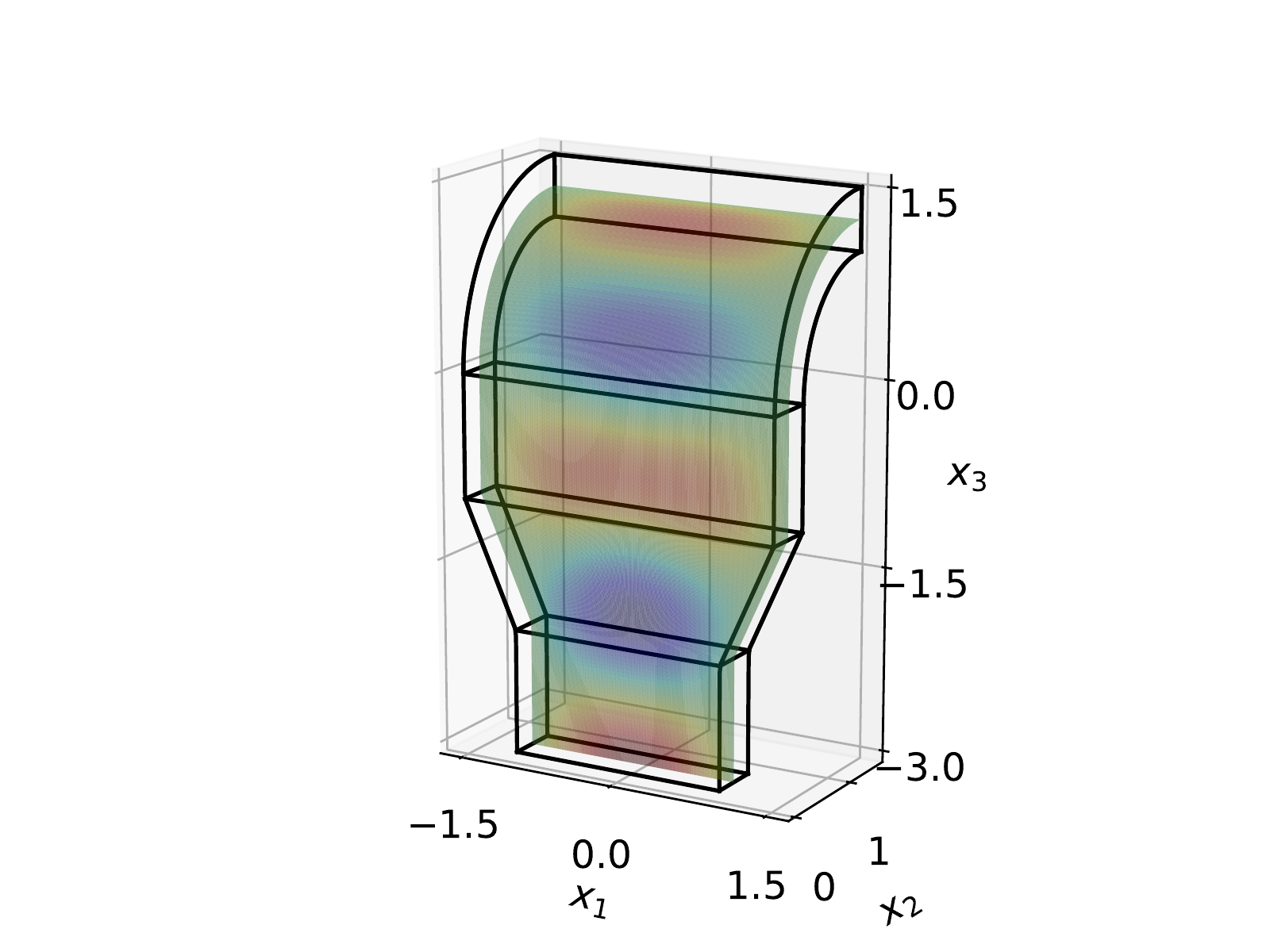}
		\caption{Solution for $\theta_1=\theta_2=\theta_3=0$ along the length and width of the waveguide for $\rho=49$.} \label{fig:wg_solution}
	\end{subfigure}
	\caption{Test case 4: Waveguide with rectangular section.}
	\label{fig:wg1}
\end{figure}

In this numerical example, the TT-IGA solver is applied to solve the scalar Helmholtz equation within a waveguide structure with parameter-dependent geometry. 
The BVP problem to be solved is
\begin{subequations}
	\begin{align}
		\Delta u(\cdot,\btheta) +\rho u(\cdot,\btheta) =0,& \quad\text{in }\: \Omega(\btheta),\\
		u(\cdot,\btheta) = g(\cdot),& \quad\text{on }\: \partial\Omega(\btheta) \cap \{(x_1,x_2,x_3):x_3=-3\},\\
		u(\cdot,\btheta) = 0 ,& \quad\text{on }\: \partial\Omega(\btheta) \setminus \{(x_1,x_2,x_3):x_3=-3\},
	\end{align}
\end{subequations}
where the parameter-dependent waveguide geometry $\Omega(\bm{\theta})$ is shown in Figure~\ref{fig:wg_parms}. Therein, the black lines represent the contour of the nominal geometry, which corresponds to the parameter values $\theta_1 = \theta_2 = \theta_3 = 0$. 
The red lines represent the maximum deformations of the geometry, corresponding to the vertices of the hypercube $\Xi=[-0.2,0.2]\times [-0.2,0.2] \times [-0.3,0.3]$.
At the boundary plane $x_3=-3$, the Dirichlet boundary condition $g(x_1,x_2)= \cos(\pi x_1) \sin(\pi x_2)$ is imposed, while the zero Dirichlet condition is assumed for the remaining boundaries. 

In Figure~\ref{fig:wg_solution}, the solution is plotted for the nominal geometry and $\rho=49$. 
For this computation, the B-spline basis consists of quadratic B-splines with dimensions $\bm{n}=(64,64,128)$. 
When resolving the parametric dependency, the parameter grid has size $\bm{\ell}=(8,8,8)$. 
For this discretization, the system matrix in the TT format is constructed within less than $3$~s. 
The TT-rank of the operator tensor is $\bm{R}=( 1 , 6, 22, 13, 14,  8 , 1 )$, corresponding to approximately $43$~MB of storage space. 
The solution has the TT-rank $\bm{R}=(1, 14, 36, 36, 36, 10, 1)$, which corresponds to less than $2$~MB of memory. 
This is less than the storage requirements for the full solution computed with an FEM solver for a fixed geometry, which exceed $4$~MB.
The TT-IGA solver is compared against the FEM solver for several geometry parameter realizations, to verify its correctness. 
Note that the GMRES solver used for the conventional FEM solution suffers from the lack of preconditioning. 
A similar increase in the computational time is also observed for the AMEn solver, however, the TT solver still remains more efficient. 
A speedup is observed if preconditioning of the local subsystems of the AMEn iterations is used.

\section{Conclusion}
\label{sec:conclusion}
This paper presented a numerical solver which uses the TT decomposition as a means to represent IGA-FEM discrete operators and solutions effectively in a low-rank tensor format for solving PDEs on parameter-dependent geometries. 
An explicit way for constructing the TT decomposition of the mass and stiffness discrete operators is provided, prior to considering any parameter dependencies.
Parameters affecting the shape of the computational domain are included in the framework by constructing a tensor product grid which combines the discretizations of the physical and the parameter space. 
To that end, the dimensionality of the solution tensor as well as of the corresponding tensor operators is extended in order to accommodate the parameter dependencies. 
Representing the combined solution tensor and tensor operators in a low-rank TT format can be very beneficial, since it can capture the dependencies between the solution and the domain deformations, an effect often referred to as the ``blessing of dimensionality'' \cite{gorban2018blessing}.
One further advantage of the TT format is the existence of the AMEn solver for handling TT-based multilinear systems and computing the solution directly in the TT format as well.
Therefore, all involved steps, i.e., constructing an IGA representation of the geometry, assembling the discrete operators, solving the multilinear system, and evaluating the solution, are performed using the TT format. 
Moreover, none of the aforementioned steps require to explicitly compute the full tensors, but merely TT approximations thereof, which are computed via TT-based cross approximation methods \cite{oseledets2010tt, savostyanov2011fast}.

The proposed TT-IGA solver has been verified using a series of numerical experiments to assess its performance and accuracy. 
The first test case consists of a convergence study considering of single parametric dependence, while the second test case examines the solver's performance for an increasing number of geometry parameters. 
With respect to the latter, the conclusion is that the framework is able to handle a moderate number of parameters, typically up to 10-12 based on our investigations. 
Two further test cases are examined, namely, a test case featuring a material discontinuity and a test case concerning waveguide simulation, where the computational advantages of the TT-IGA solver are showcased against classical IGA or FEM-based solution methods.
In both test cases, the AMEn solver is faster for finer discretizations, even if compared to traditional iterative solvers that resolve a single geometry realization.
Regarding the storage requirements, the TT format is found to be particularly efficient, as it is able to store the parameter-dependent solution with less entries than the full tensor needed for a single geometry configuration by a standard solver.
This result showcases the expressive power of the TT format for representing multidimensional structured data.

Future work in this direction will consider the addition of NURBS-based parametrizations with parameter-dependent knots and weights, as the present work considered B-spline basis functions only.
Additionally desirable would be an extension of the method suggested in \cite{pan2019low}, such that the construction of IGA representations of computational domains in the low-rank TT format starting from boundary patches is applicable for the case of parameter-dependent domains.

\subsubsection*{Software}
\noindent The code containing the implementation and the numerical tests is publicly available at: \\ 
\texttt{https://https://github.com/ion-g-ion/code-paper-tt-iga}.

\subsubsection*{Aknowledgement}
\noindent All authors are supported by the Graduate School
Computational Engineering within the Centre for Computational Engineering at the Technische Universit\"at Darmstadt. 

\bibliography{refs}{}
\bibliographystyle{plain}
\end{document}